\documentclass[11pt]{article}%
\usepackage{amsmath}
\usepackage{amsthm}
\usepackage{dsfont}
\usepackage{graphicx}
\usepackage{mathdots}
\usepackage{xcolor}%
\usepackage{amsfonts}%
\usepackage{amssymb}
\newcommand{\eqnb}{\begin{equation}}
\newcommand{\eqne}{\end{equation}}

\newtheorem{The}{Theorem}

\newtheorem{Lem}{Lemma}

\newtheorem{Rem}{Remark}

\begin{document}

\title{\textbf{Matched Queues with Matching Batch Pair $(m, n)$}}
\author{Heng-Li Liu$^{a}$, Quan-Lin Li$^{b}$ \thanks{Corresponding author: Q.L. Li
(liquanlin@tsinghua.edu.cn)}, Chi Zhang$^{b}$,\\$^{a}$School of Economics and Management Sciences \\Yanshan University, Qinhuangdao 066004, China\\$^{b}$School of Economics and Management \\Beijing University of Technology, Beijing 100124, China}
\maketitle

\begin{abstract}
In this paper, we develop the matrix-analytic method to discuss an interesting
but challenging bilateral stochastic matching problem: A matched queue with
matching batch pair $(m, n)$ and two types of impatient customers, where the
two types of customers arrive according to two independent Poisson processes.
Once $m$ A-customers and $n$ B-customers are matched as a group, the $m + n$
customers immediately leave the system. We show that this matched queue can be
expressed as a novel bidirectional level-dependent quasi-birth-and-death (QBD)
process whose analysis has its own interests, and specifically, computing the
maximal non-positive inverse matrices of bidirectional infinite sizes by using
the RG-factorizations. Based on this, we can provide an effective
matrix-analytic method to deal with this matched queue, including the system
stability, the average stationary queue lengths, the average sojourn times,
and the departure process. We believe that the methodology and results
developed in this paper can be applicable to studying more general matched
queueing systems, which are widely encountered in many practical areas.


\end{abstract}

\section{Introduction}

In this paper, we consider an interesting but challenging matched queue with
matching batch pair $(m,n)$ and two types (i.e., types A and B) of impatient
customers, where $m$ A-customers and $n$ B-customers are matched as a group
and the $m+n$ customers leave the system immediately. To our best knowledge,
this paper is the first to study the more general matched queue due to the
matching batch pair $(m,n)$. It is worthwhile to note that the matching batch
pair $(m,n)$\ makes not only the model more suitable for many practically
matching needs but also our Markov modeling and analysis more challenging. We
show that this matched queue can be expressed as a level-dependent QBD process
with bidirectional infinite levels, and thus develop some new theory of
level-dependent QBD processes with bidirectional infinite levels, such as the
system stability, the stationary probability vectors, the sojourn times, the
first passage times, and the departure processes through dealing with the
matrices of bidirectional infinite sizes by using the RG-factorizations given
in Li \cite{Li:2010}. Based on this, we can provide a detailed analysis for
this matched queue, including the system stability, the average stationary
queue lengths, the average sojourn times, and the departure process. Also, we
can further develop some effective algorithms for analyzing the performance
measures of this matched system.

So far, more and more matching problems (e.g., the double-ended queues, the
matched queues, and more generally, the fork-join queues) have been widely
encountered in many different practical areas, for example, sharing economy,
ridesharing platform, bilateral market, organ transplantation, assembly
systems, taxi services, seaport and airport, and so on. Important research
examples include: \textit{Organ transplantation} by Zenios \cite{Zen:1999},
Boxma et al. \cite{Box:2011}, Stanford et al. \cite{Sta:2014}, and Elalouf et
al. \cite{Ela:2018}; \textit{taxi services} by Giveen \cite{Giv:1961,
Giv:1963}, Kashyap \cite{Kas:1965, Kas:1966, Kas:1967}, Bhat \cite{Bah:1970},
Baik et al. \cite{Bai:2002}, Shi and Lian \cite{Shi:2016}, and Zhang et al.
\cite{Zha:2019}; \textit{baggage claim} by Browne et al. \cite{Bro:1970};
\textit{sharing economy} by Cheng \cite{Cheng:2016}, Sutherland and Jarrahi
\cite{Sut:2018}, and Benjaafar and Hu \cite{Ben:2020}; \textit{assembly
systems} by Hopp and Simon \cite{Hop:1989}, Som et al. \cite{Som:1994}, and
Ramachandran and Delen \cite{Ram:2005}; \textit{health care} by Pandey and
Gangeshwer \cite{Pan:2018}; \textit{multimedia synchronization} by Steinmetz
\cite{Ste:1990} and Parthasarathy et al. \cite{Par:1999}; and so forth.

Besides these, in recent years, an emerging hot research topic of matched
queues has been focused on \textit{ridesharing platform}. Many ridesharing
companies spring up as a result of rapid development of mobile networks, smart
phones and location technologies, for example, Uber in transportation, Airbnb
in housing, Eatwith in eating, Rent the Runway in dressing, and so on. Readers
may refer to, such as Azevedo and Weyl \cite{Aze:2016}, Duenyas et al.
\cite{Due:1997}, Hu and Zhou \cite{Hu:2015}, Banerjee and Johari
\cite{Ban:2019} and Braverman et al. \cite{Bra:2019}. Obviously, the matched
queues have become useful and necessary in both theory research and real
applications of ridesharing platforms with various different services.

Now, we summarize the literature of matched queues from three different
aspects of matching batch pair $(m,n)$ as follows:

\textbf{The matching batch pair }$(1,1)$\textbf{.} Early research of matching
queues first focused on some simple double-ended (or matched) systems with
matching batch pair $(1,1)$. Also, the matched queues with matching batch pair
$(1,1)$ have attracted numerous researchers' attention since a pioneering work
by Kendall \cite{Ken:1951}, and crucially, some effective methodologies and
available results have been developed from multiple research perspectives
listed below.

\textit{The Markov process:} Such a process was the first effective method
employed in early research of matched queues. For a simple matched queue,
Sasieni \cite{Sas:1961}, Giveen \cite{Giv:1961} and Dobbie \cite{Dob:1961}
established the Chapman-Kolmogorov forward differential-difference equations,
in which the customers' impatient behavior was introduced to guarantee the
system stability. Since then, the Markov process analysis of matched queues
was further developed from two different research lines:

(a) A finite state space. When the two waiting rooms of the matched queue are
both finite, the Markov process is established on a finite state space. In
this case, Jain \cite{Jai:1962} and Kashyap \cite{Kas:1965, Kas:1966,
Kas:1967} applied the supplementary variable method to be able to deal with
the matched queue with a Poisson arrival process and a renewal arrival
process, and specifically, Takahashi et al. \cite{Tak:2000} considered the
matched queue with a Poisson arrival process and a PH-renewal arrival process.
In addition, Sharma and Nair\cite{Sha:1991} used the matrix theory to analyze
the transient behavior of a Markovian matched queue. Chai et al.
\cite{Chai:2019} considered a batch matching queueing system with impatient
servers and bounded rational customers, and each server serves the customers
in batches with finite service capacity.

(b) An infinite state space. When the two waiting rooms of the matched queue
are both infinite, the Markov process is set up with a bidirectional infinite
state space. In general, it is always difficult to analyze such a Markov
process on a bidirectional infinite state space. Latouche \cite{Lat:1981}
applied the matrix-geometric solution to analyze several bilateral matched
queues with paired input. Conolly et al. \cite{Con:2002} applied the Laplace
transform to discuss the time-dependent performance measures of the matched
queue with state-dependent impatience. Di Crescenzo et al. \cite{Di:2012,
Di:2018} discussed the transient and stationary probabilities of a
time-nonhomogeneous matched queue with catastrophes and repairs. Diamant and
Baron \cite{Dia:2019} analyzed a matched queue with priority and impatient customers.

\textit{The fluid and diffusion approximations:} In a matched queue, if the
arrivals of A- and B-customers are both general renewal processes, then the
fluid and diffusion approximations become an effective (but approximative)
method. Jain \cite{Jai:2000} applied the diffusion approximation to discuss
the G$^{\text{X}}$/G$^{\text{Y}}$/1 matched queue. Di Crescenzo et al.
\cite{Di:2012, Di:2018} discussed a matched queue by means of a jump-diffusion
approximation. Liu et al. \cite{Liu:2014} discussed some diffusion models for
the matched queues with renewal arrival processes. B\"{u}ke and Chen
\cite{Buke:2017} applied the fluid and diffusion approximations to study the
probabilistic matching systems. Liu \cite{Liu:2019} used the diffusion
approximation to analyze the matched queues with reneging in heavy traffic.

\textit{Other effective methods:} Adan et al. \cite{Adan:2018a, Adan:2018b}
and Visschers et al. \cite{Vis:2019} discussed the matched systems with
multi-type jobs and multi-type servers by using the product solution of
queueing networks. Kim et al. \cite{Kim:2010} provided a simulation model to
analyze a more general matched queue. Jain \cite{Jai:1995} proposed a sample
path analysis for studying the matched queue with time-dependent rates.
Af\`{e}che et al. \cite{Afe:2014} applied the level-cross method to discuss
the batch matched queue with abandonment. Wu and He \cite{Wu:2020} used the
multi-layer Markov modulated fluid flow (MMFF) processes to deal with a
double-sided queueing model with marked Markovian arrival processes and finite
discrete abandonment times.

\textit{Control of matched queues: }Hlynka and Sheahan \cite{Hly:1987}
analyzed the control rates in a matched queue with two Poisson inputs. Gurvich
and Ward \cite{Gur:2014} discussed dynamic control of the matching queues.
B\"{u}ke and Chen \cite{Buke:2015} studied stabilizing admission control
policies for the probabilistic matching systems. Lee et al. \cite{Lee:2019}
studied optimal control of a time-varying double-ended production queueing system.

\textbf{The matching batch pair }$(1,n)$\textbf{.} As a key generalization, Xu
et al. \cite{Xu:1990} first discussed a matched queue with matching batch pair
$(1,n)$, in which for the two waiting rooms, one is finite while another is
infinite. Under two Poisson inputs and a PH service time distribution, they
applied the matrix-geometric solution to obtain the stability condition of the
system, and to study the stationary queue lengths for the both classes of
customers. Since then, further research includes Xu and He \cite{Xu:1993a,
Xu:1993b}. Yuan \cite{Yuan:1992} applied Markov chains of M/G/1 type to
consider a matched queue with matching batch pair $(1,n)$ and under two
Poisson inputs and a general service time distribution. Li and Cao
\cite{Li:1996} discussed the matched queue with matching batch pair $(1,n)$
and under two batch Markovian arrival processes (BMAPs) and a general service
time distribution.

\textbf{The matching batch pair }$(m,n)$\textbf{.} To our best knowledge, this
paper is the first to study the matched queues with matching batch pair
$(m,n)$, where the two waiting rooms are both infinite.\ We express this
matched queue as a level-dependent QBD process with bidirectional infinite
levels, and apply the RG-factorizations given in Li \cite{Li:2010} to obtain
the average stationary queue lengths, the average sojourn times, and the
departure process. Note that Liu et al. \cite{Liu:2020} is a closely related
work to study such matched queues whose corresponding Markov processes are
block-structured and level-dependent. Different from Liu et al.
\cite{Liu:2020}, the matching batch pair $(m,n)$ makes the Markov block
structure of bidirectional infinity sizes more challenging. In addition, we
develop the matrix-analytic method to study the average sojourn time and the
departure process through setting up a new PH distribution of bidirectional
infinity sizes and a new Markovian arrival process with marked transitions
(MMAP) of bidirectional infinity sizes, respectively.

In what follows, it is necessary to discuss some random features of the
matched queues from the theory of Markov processes.

On the one hand, the matched queues are a type of interesting and classic
systems in early research of queuing systems, but their available
methodologies and results are fewer than those developed for other types of
classic queueing systems, for example, processor-sharing queues (Yashkov
\cite{Yas:1987} and Yashkov and Yashkov \cite{Yas:2007}), and retrial queues
(Falin and Templeton \cite{Fal:1997} and Artalejo and G\'{o}mez-Corral
\cite{Art:2008}). In fact, the matrix-analytic method (based on the
level-dependent Markov processes) have been applied to analysis of the retrial
queues or the processor-sharing queues, e.g., See Chapter 5 in Artalejo and
G\'{o}mez-Corral \cite{Art:2008} and Chapter 7 in Li \cite{Li:2010}. On the
other hand, it is relatively difficult to discuss the queues with impatient
customers (even though the impatient times are exponential). Note that the
general impatient times can make the embedded Markov process analysis of the
queue very complicated and even impossible except for using the fluid and
diffusion approximations under an approximate goal, e.g., see Boots and Tijms
\cite{Boo:1999}, Zeltyn and Mandelbaum \cite{Zel:2005}, and Puha and Ward
\cite{Puh:2019}. Furthermore, the customers' impatient behavior greatly
complicates the analysis of matched queues due to the level-dependent
structure of their corresponding Markov process with bidirectional infinite
levels. Thus, the matrix-analytic method needs to further be developed through
applying the Markov process with bidirectional infinite levels to dealing with
the matched queueing example with matching batch pair $(m,n)$.

In the study of matched queues, the matching batch pair $(m,n)$ greatly
complicates how to concretely write the infinitesimal generator of the
level-dependent QBD process with bidirectional infinite levels. Also, the
matching batch pair $(m,n)$ convincingly motivate us to develop the
matrix-analytic method of matched queues, which can be applied to dealing with
many practical matching problems. For the matrix-analytic method to the study
of matched queues, the level-dependent Markov processes are a key and their
analysis is based on the RG-factorizations given in Li \cite{Li:2010}. Thus,
the RG-factorizations play a key role in the study of matched queues. By using
the RG-factorizations, we can further develop some effective algorithms (also
see some algorithmic research by Bright and Taylor \cite{Bri:1995, Bri:1997},
Takine \cite{Tak:2016}\ and Liu et al. \cite{Liu:2020}) to be able to
numerically analyze performance measures of the matched queues. Although some
results given in this paper are regarded as superficially coming from Chapter
2 of Li \cite{Li:2010}, the level-dependent block structure leads to some new
advances in the matrix-analytic method of matched queues, including the
level-dependent Markov processes with bidirectional infinite levels, the
maximal non-positive inverse matrices of bidirectional infinite sizes, the
sojourn times of bidirectional infinite sizes, the departure processes of
bidirectional infinite sizes, and so forth.

The Markovian arrival process (MAP) is a useful mathematical tool, for
example, for describing bursty traffic and dependent arrivals in many real
systems, such as computer and communication networks, manufacturing systems,
transportation networks and so on. Readers may refer to, such as Chapter 5 in
Neuts \cite{Neu:1989}, Lucantoni \cite{Luc:1991}, Chapter 1 in \cite{Li:2010}
and references therein. Further, He \cite{He:1996} and He and Neuts
\cite{He:1998} introduced the MMAP, which is useful in modeling input (or
departure) processes of stochastic systems with several types of items (e.g.,
customers or orders). In this paper, we use the MMAP of bidirectional infinity
sizes to study the departure process with three types of customers in the
matched queue with matching batch pair $(m,n)$.

Based on the above analysis, we summarize the main contributions of this paper
as follows:

\begin{itemize}
\item[(1)] We describe and analyze a more general matched queue with matching
batch pair $(m,n)$, where $m$ A-customers and $n$ B-customers are matched as a
group and the $\left(  m+n\right)  $ customers leave the system immediately.
Such a matched queue can widely be used to study many practically matching
problems, for example, sharing economy, ridesharing platform, bilateral
market, organ transplantation, taxi services, assembly systems, and so on.

\item[(2)] We express the matched queue with matching batch pair $(m,n)$ as a
level-dependent QBD process with bidirectional infinite levels. Thus, we
develop the matrix-analytic method by means of the level-dependent QBD
processes with bidirectional infinite levels, such as the system stability,
the stationary probability vectors, the sojourn times, the first passage
times, and the departure processes. A key of our method is based on applying
the RG-factorizations given in Li \cite{Li:2010} to dealing with the maximal
non-positive inverse matrices of bidirectional infinite sizes.

\item[(3)] We provide a detailed analysis for the matched queue with matching
batch pair $(m,n)$, including the system stability, the average stationary
queue lengths, the average sojourn times, and the departure process.
Specifically, the average sojourn times are given a better upper bound by
using a new phase-type distribution of bidirectional infinity sizes, and the
departure process with three types of customers is established in terms of the
MMAP of bidirectional infinity sizes. Also, some numerical examples are used
to indicate our theoretical results.
\end{itemize}

\vskip
0.4cm

The structure of this paper is organized as follows. Section 2 describes a
more general matched queue with matching batch pair $(m,n)$ and two types of
impatient customers. Section 3 expresses this matched queue as a
level-dependent QBD process with bidirectional infinite levels, and provides
the stability condition of the system. Section 4 studies the stationary
probability vector of the QBD process with bidirectional infinite levels, and
thus computes the average stationary queue length of any A- or B-customer.
Section 5 computes the average sojourn time of any A- or B-customer by using
three different techniques: The Little's formula, a probabilistic calculation,
and an upper bound, respectively. Section 6 uses the MMAP of bidirectional
infinite sizes to discuss the departure process with three types of customers.
Finally, some concluding remarks are given in Section 7.

\section{Model Description}

In this section, we describe a more general matched queue with matching batch
pair $(m,n)$ and two types of impatient customers, and also introduce
operational mechanism, system parameters and basic notation.

In the matched queue with matching batch pair $(m,n)$, $m$ A-customers and $n$
B-customers are matched as a group which leaves the system immediately once
such a matching is successful, and also the customers' impatient behavior is
used to guarantee the stability of the system. Figure 1 provides a physical
illustration for such a matched queue.

\begin{figure}[tbh]
\centering                      \includegraphics[width=14cm]{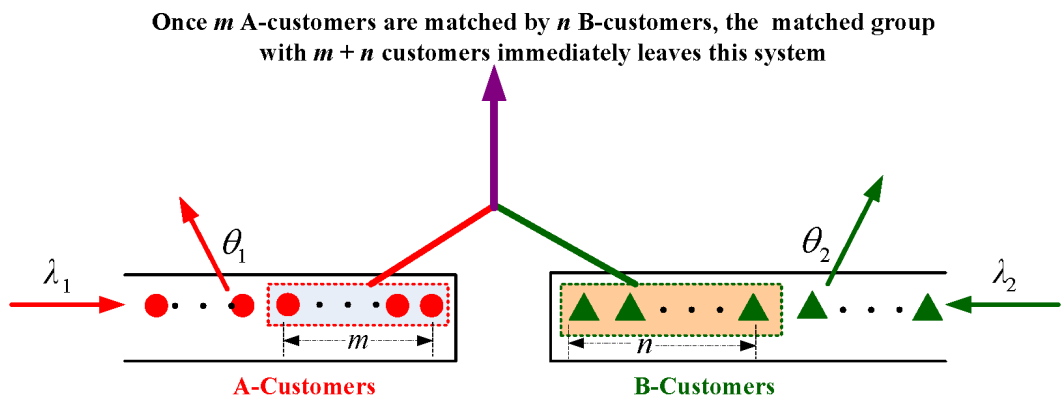}
\caption{A physical illustration of the matched queue}%
\end{figure}

Now, we provide a more detailed description for the matched queue as follows:

\textbf{(1) Arrival processes.} The A- and B-customers arrive at the queueing
system according to two Poisson processes with rates $\lambda_{1}$ and
$\lambda_{2}$, respectively.

\textbf{(2) Matching processes.} When $m$ customers of type A and $n$
customers of type B are present, there is a match and these $m+n$ customers
immiediately leave the system for $m,n\geq1$. Customers for whom there is not
yet a match must wait. Their matching process follows a First-Come-First-Match
discipline. We assume that the two waiting spaces of A- and B-customers are
both infinite.

\textbf{(3) Impatient behavior.} If an A-customer (resp. a B-customer) stays
at the queueing system for a long time, then she will have some impatient
behavior. We assume that the impatient time of an A-customer (resp. a
B-customer) is exponentially distributed with impatient rate $\theta_{1}$
(resp. $\theta_{2}$) for $\theta_{1},\theta_{2}>0$.

We assume that all the random variables defined above are independent of each other.

\begin{Rem}
(a) The customers' impatient behavior given in Assumption (3) is used to
ensure the stability of the matched queue.

(b) The matching discipline given in Assumption (2) indicates that more than
$m$ A-customers and more than $n$ B-customers cannot simultaneously exist in
their waiting spaces.
\end{Rem}

\section{A QBD Process with Bidirectional Infinite Levels}

In this section, we describe the matched queue with matching batch pair
$(m,n)$ as a new bidirectional level-dependent QBD process to express , and
obtain a sufficient condition under which this matched queue is stable.

We denote by $N_{1}\left(  t\right)  $ and $N_{2}\left(  t\right)  $ the
numbers of A- and B-customers in the matched queue at time $t\geq0$,
respectively. Then the matched queue with matching batch pair $(m,n)$ is
related to a two-dimensional Markov process $\left\{  \left(  N_{1}\left(
t\right)  ,N_{2}\left(  t\right)  \right)  ,\text{ }t\geq0\right\}  $. Note
that once $m$ A-customers and $n$ B-customers are matched as a group, the
$m+n$ customers immediately leave the queueing system, thus more than $m$
A-customers and more than $n$ B-customers cannot simultaneously exist in their
waiting spaces. Based on this, the state space of the Markov process $\left\{
\left(  N_{1}\left(  t\right)  ,N_{2}\left(  t\right)  \right)  ,\text{ }%
t\geq0\right\}  $ is given by%
\begin{align*}
\Omega=  &  \left\{  \left(  i,j\right)  :0\leq i\leq m-1,0\leq j\leq
n-1\right\}  \cup\left\{  \left(  i,j\right)  :i\geq m,0\leq j\leq n-1\right\}
\\
&  \cup\left\{  \left(  i,j\right)  :0\leq i\leq m-1,j\geq n\right\}  .
\end{align*}
In this case, we write%
\begin{align}
\text{Level }0=  &  \left\{  \left(  0,0\right)  ,\left(  0,1\right)
,\ldots,\left(  0,n-1\right)  ;\left(  1,0\right)  ,\text{ }\left(
1,1\right)  ,,\ldots,\text{\ }\left(  1,n-1\right)  ;\right. \label{L0}\\
&  \left.  \ldots;\text{\ }\left(  m-1,0\right)  ,\text{\ }\left(
m-1,1\right)  ,\ldots,\text{\ }\left(  m-1,n-1\right)  \right\}  ,\nonumber
\end{align}
for $k\geq1,$
\begin{align}
\text{Level }k=  &  \left\{  \left(  km,0\right)  ,\left(  km,1\right)
,\ldots,\left(  km,n-1\right)  ;\left(  km+1,0\right)  ,\text{ }\left(
km+1,1\right)  ,\ldots,\left(  km+1,n-1\right)  ;\right. \label{Lk}\\
&  \left.  \ldots;\text{\ }\left(  km+\left(  m-1\right)  ,0\right)  ,\left(
km+\left(  m-1\right)  ,1\right)  ,\ldots,\left(  km+\left(  m-1\right)
,n-1\right)  \right\}  ,\nonumber
\end{align}
and $l\leq-1,$%
\begin{align}
\text{Level }l=  &  \left\{  \left(  m-1,\left(  -l\right)  n+n-1\right)
,\left(  m-2,\left(  -l\right)  n+n-1\right)  ,\ldots,\left(  0,\left(
-l\right)  n+n-1\right)  ;\ldots;\right. \nonumber\\
&  \left(  m-1,\left(  -l\right)  n+1\right)  ,\left(  m-2,\left(  -l\right)
n+1\right)  ,\ldots,\left(  0,\left(  -l\right)  n+1\right)  ;\nonumber\\
&  \left.  \left(  m-1,\left(  -l\right)  n\right)  ,\left(  m-2,\left(
-l\right)  n\right)  ,\ldots,\left(  0,\left(  -l\right)  n\right)  \right\}
, \label{Ll}%
\end{align}
Therefore, we have%
\[
\Omega=\bigcup\limits_{k=-\infty}^{\infty}\text{Level }k.
\]

\textbf{Example one:} As an illustrated example, we take $m=2$ and $n=3$. In
this case, the state transition relations of Markov process $\left\{  \left(
N_{1}\left(  t\right)  ,N_{2}\left(  t\right)  \right)  ,\text{ }%
t\geq0\right\}  $ are depicted in Figure 2. Also, we observe that each level
is a state set formed by many states in a rectangle (i.e., multiple state lines).

\begin{figure}[tbh]
\centering                      \includegraphics[width=14cm]{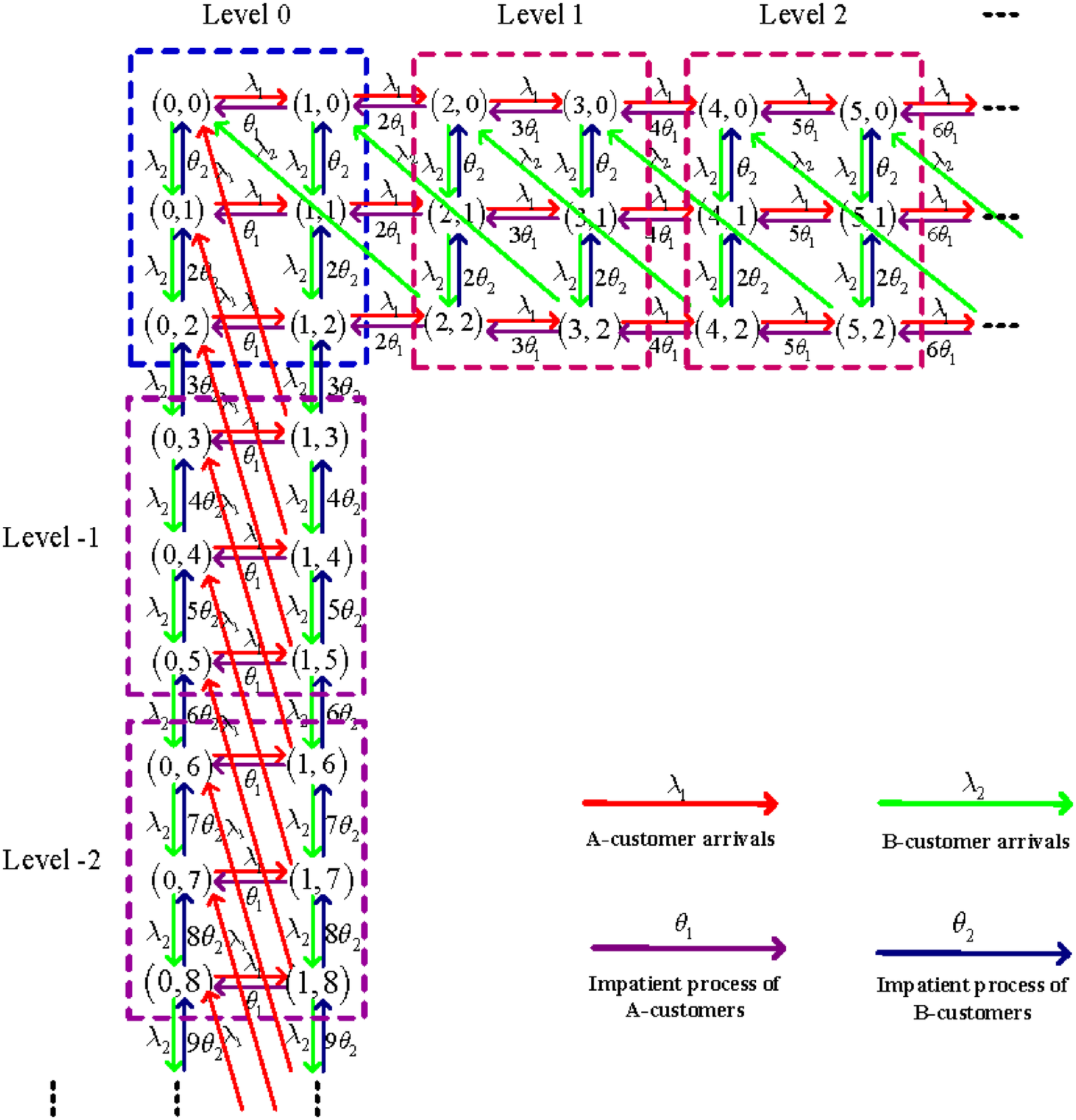}
\caption{The state transition relations of the bilateral QBD process}%
\end{figure}

From Levels $k$ for $-\infty<k<\infty$ or Figure 2, it is easy to see that the
Markov process $\{(N_{1}\left(  t\right)  ,N_{2}\left(  t\right)  ),$
$t\geq0\}$ is a new level-dependent QBD process with bidirectional infinite
levels whose infinitesimal generator is given by
\begin{equation}
Q=\left(
\begin{array}
[c]{ccccccccccc}%
\ddots & \ddots & \ddots &  &  &  &  &  &  &  & \\
& B_{0}^{\left(  -3\right)  } & B_{1}^{\left(  -3\right)  } & B_{2}^{\left(
-3\right)  } &  &  &  &  &  &  & \\
&  & B_{0}^{\left(  -2\right)  } & B_{1}^{\left(  -2\right)  } &
B_{2}^{\left(  -2\right)  } &  &  &  &  &  & \\
&  &  & B_{0}^{\left(  -1\right)  } & B_{1}^{\left(  -1\right)  } &
\fbox{$B_{2}^{\left(  -1\right)  }$} &  &  &  &  & \\
&  &  &  & \fbox{$B_{0}^{\left(  0\right)  }$} & C & \fbox{$A_{0}^{\left(
0\right)  }$} &  &  &  & \\
&  &  &  &  & \fbox{$A_{2}^{\left(  1\right)  }$} & A_{1}^{\left(  1\right)  }
& A_{0}^{\left(  1\right)  } &  &  & \\
&  &  &  &  &  & A_{2}^{\left(  2\right)  } & A_{1}^{\left(  2\right)  } &
A_{0}^{\left(  2\right)  } &  & \\
&  &  &  &  &  &  & A_{2}^{\left(  3\right)  } & A_{1}^{\left(  3\right)  } &
A_{0}^{\left(  3\right)  } & \\
&  &  &  &  &  &  &  & \ddots & \ddots & \ddots
\end{array}
\right)  , \label{Q}%
\end{equation}
where, for $k\geq0,$
\begin{align*}
\text{Level }k=  &  \left\{  \left(  km,0\right)  ,\left(  km,1\right)
,\ldots,\left(  km,n-1\right)  ;\left(  km+1,0\right)  ,\text{ }\left(
km+1,1\right)  ,\ldots,\left(  km+1,n-1\right)  ;\right. \\
&  \left.  \ldots;\text{\ }\left(  km+\left(  m-1\right)  ,0\right)  ,\left(
km+\left(  m-1\right)  ,1\right)  ,\ldots,\left(  km+\left(  m-1\right)
,n-1\right)  \right\}  ,
\end{align*}
we have
\[
C=\left(
\begin{array}
[c]{ccccc}%
A_{1,1}^{\left(  0\right)  } & A_{2,1} &  &  & \\
A_{3,1}^{\left(  0\right)  } & A_{1,2}^{\left(  0\right)  } & A_{2,2} &  & \\
& \ddots & \ddots & \ddots & \\
&  & A_{3,m-2}^{\left(  0\right)  } & A_{1,m-1}^{\left(  0\right)  } &
A_{2,m-1}\\
&  &  & A_{3,m-1}^{\left(  0\right)  } & A_{1,m}^{\left(  0\right)  }%
\end{array}
\right)  _{mn\times mn},\text{ \ }%
\]%
\[
A_{1,i}^{\left(  0\right)  }=\left(
\begin{array}
[c]{ccccc}%
a_{1,i}^{\left(  0,1\right)  } & \lambda_{2} &  &  & \\
\theta_{2} & a_{1,i}^{\left(  0,2\right)  } & \lambda_{2} &  & \\
& \ddots & \ddots & \ddots & \\
&  & \left(  n-2\right)  \theta_{2} & a_{1,i}^{\left(  0,n-1\right)  } &
\lambda_{2}\\
&  &  & \left(  n-1\right)  \theta_{2} & a_{1,i}^{\left(  0,n\right)  }%
\end{array}
\right)  _{n\times n},\text{ }1\leq i\leq m,
\]%
\[
a_{1,i}^{\left(  0,r\right)  }=-\left(  \lambda_{1}+\lambda_{2}+\left(
r-1\right)  \theta_{2}+\left(  i-1\right)  \theta_{1}\right)  ,\text{ }1\leq
r\leq n;
\]%
\[
A_{1}^{\left(  k\right)  }=\left(
\begin{array}
[c]{ccccc}%
A_{1,1}^{\left(  k\right)  } & A_{2,1} &  &  & \\
A_{3,1}^{\left(  k\right)  } & A_{1,2}^{\left(  k\right)  } & A_{2,2} &  & \\
& \ddots & \ddots & \ddots & \\
&  & A_{3,m-2}^{\left(  k\right)  } & A_{1,m-1}^{\left(  k\right)  } &
A_{2,m-1}\\
&  &  & A_{3,m-1}^{\left(  k\right)  } & A_{1,m}^{\left(  k\right)  }%
\end{array}
\right)  _{mn\times mn},\text{ \ }k\geq0,
\]%
\[
A_{1,i}^{\left(  k\right)  }=\left(
\begin{array}
[c]{ccccc}%
a_{1,i}^{\left(  k,1\right)  } & \lambda_{2} &  &  & \\
\theta_{2} & a_{1,i}^{\left(  k,2\right)  } & \lambda_{2} &  & \\
& \ddots & \ddots & \ddots & \\
&  & \left(  n-2\right)  \theta_{2} & a_{1,i}^{\left(  k,n-1\right)  } &
\lambda_{2}\\
&  &  & \left(  n-1\right)  \theta_{2} & a_{1,i}^{\left(  k,n\right)  }%
\end{array}
\right)  _{n\times n},\text{ }1\leq i\leq m,
\]%
\[
a_{1,i}^{\left(  k,r\right)  }=-\left(  \lambda_{1}+\lambda_{2}+\left(
r-1\right)  \theta_{2}+\left(  km+i-1\right)  \theta_{1}\right)  ,\text{
}1\leq r\leq n,
\]%
\[
A_{2,i}=\left(
\begin{array}
[c]{cccc}%
\lambda_{1} &  &  & \\
& \lambda_{1} &  & \\
&  & \ddots & \\
&  &  & \lambda_{1}%
\end{array}
\right)  _{n\times n},1\leq i\leq m-1,
\]%
\[
A_{3,i}^{\left(  k\right)  }=\left(
\begin{array}
[c]{cccc}%
\left(  km+i\right)  \theta_{1} &  &  & \\
& \left(  km+i\right)  \theta_{1} &  & \\
&  & \ddots & \\
&  &  & \left(  km+i\right)  \theta_{1}%
\end{array}
\right)  _{n\times n},\text{ }1\leq i\leq m-1;
\]
note that $A_{0}^{\left(  k\right)  }$ is a transition rate matrix from Level
$k$ to Level $k+1$, we obtain
\[
A_{0}^{\left(  k\right)  }=\left(
\begin{array}
[c]{cc}%
{\huge 0}_{\left[  \left(  m-1\right)  n\right]  \times n} & {\huge 0}%
_{\left[  \left(  m-1\right)  n\right]  \times\left[  \left(  m-1\right)
n\right]  }\\
\text{ }\left.
\begin{array}
[c]{ccc}%
\lambda_{1} &  & \\
& \ddots & \\
&  & \lambda_{1}%
\end{array}
\right\}  n\times n & {\huge 0}_{n\times\left[  \left(  m-1\right)  n\right]
}%
\end{array}
\right)  _{mn\times mn},\text{\ }k\geq0;
\]
since $A_{2}^{\left(  k\right)  }$ is a transition rate matrix from Level $k$
to Level $k-1$, we obtain%
\[
A_{2}^{\left(  k\right)  }=\left(
\begin{array}
[c]{ccc}%
\left.
\begin{array}
[c]{ccc}%
0 &  & \\
& \ddots & \\
\lambda_{2} &  & 0
\end{array}
\right\}  n\times n & \cdots & \left.
\begin{array}
[c]{ccc}%
km\theta_{1} &  & \\
& \ddots & \\
&  & km\theta_{1}%
\end{array}
\right\}  n\times n\\
& \ddots & \vdots\\
&  & \left.
\begin{array}
[c]{ccc}%
\text{ \ }0\text{ \ } & \text{ \ \ \ } & \text{ \ \ \ \ \ \ \ }\\
\text{ \ } & \ddots & \\
\lambda_{2} &  & 0
\end{array}
\right\}  n\times n
\end{array}
\right)  _{mn\times mn},\text{\ }k\geq1.
\]
For $l\leq-1,$%
\begin{align*}
\text{Level }l=  &  \left\{  \left(  m-1,\left(  -l\right)  n+n-1\right)
,\left(  m-2,\left(  -l\right)  n+n-1\right)  ,\ldots,\left(  0,\left(
-l\right)  n+n-1\right)  ;\ldots;\right. \\
&  \left(  m-1,\left(  -l\right)  n+1\right)  ,\left(  m-2,\left(  -l\right)
n+1\right)  ,\ldots,\left(  0,\left(  -l\right)  n+1\right)  ;\\
&  \left.  \left(  m-1,\left(  -l\right)  n\right)  ,\left(  m-2,\left(
-l\right)  n\right)  ,\ldots,\left(  0,\left(  -l\right)  n\right)  \right\}
,
\end{align*}
we have%
\[
B_{1}^{\left(  l\right)  }=\left(
\begin{array}
[c]{ccccc}%
B_{1,1}^{\left(  l\right)  } & B_{2,1}^{\left(  l\right)  } &  &  & \\
B_{3,1} & B_{1,2}^{\left(  l\right)  } & B_{2,2}^{\left(  l\right)  } &  & \\
& \ddots & \ddots & \ddots & \\
&  & B_{3,n-2} & B_{1,n-1}^{\left(  l\right)  } & B_{2,n-1}^{\left(  l\right)
}\\
&  &  & B_{3,n-1} & B_{1,n}^{\left(  l\right)  }%
\end{array}
\right)  _{mn\times mn},\text{ }l\leq-1,
\]%
\[
B_{1,i}^{\left(  l\right)  }=\left(
\begin{array}
[c]{ccccc}%
b_{1,i}^{\left(  l,1\right)  } & \left(  m-1\right)  \theta_{1} &  &  & \\
\lambda_{1} & b_{1,i}^{\left(  l,2\right)  } & \left(  m-2\right)  \theta_{1}
&  & \\
& \ddots & \ddots & \ddots & \\
&  & \lambda_{1} & b_{1,i}^{\left(  l,m-1\right)  } & \theta_{1}\\
&  &  & \lambda_{1} & b_{1,i}^{\left(  l,m\right)  }%
\end{array}
\right)  _{m\times m},\text{ }1\leq i\leq n,
\]%
\[
b_{1,i}^{\left(  l,r\right)  }=-\left(  \lambda_{1}+\lambda_{2}+\left(
m-r\right)  \theta_{1}+\left(  -ln+n-i\right)  \theta_{2}\right)  ,\text{
}1\leq r\leq m,
\]%
\[
B_{2,i}^{\left(  l\right)  }=\left(
\begin{array}
[c]{cccc}%
\left(  -ln+n-i\right)  \theta_{2} &  &  & \\
& \left(  -ln+n-i\right)  \theta_{2} &  & \\
&  & \ddots & \\
&  &  & \left(  -ln+n-i\right)  \theta_{2}%
\end{array}
\right)  _{m\times m},1\leq i\leq n-1,
\]%
\[
B_{3,i}=\left(
\begin{array}
[c]{cccc}%
\lambda_{2} &  &  & \\
& \lambda_{2} &  & \\
&  & \ddots & \\
&  &  & \lambda_{2}%
\end{array}
\right)  _{m\times m},1\leq i\leq n-1;
\]
note that $B_{0}^{\left(  0\right)  }$ is a transition rate matrix from Level
$0$ to Level $-1$ while Level $0$ and Level $-1$ have a different
lexicographic order, we obtain%

\[
B_{0}^{\left(  0\right)  }=\left(
\begin{array}
[c]{cc}%
{\huge 0}_{n\times\left[  \left(  n-1\right)  m\right]  } & \left.
\begin{array}
[c]{ccc}%
0 & \cdots & 0\\
\vdots & \ddots & \vdots\\
0 & \cdots & \lambda_{2}%
\end{array}
\right\}  n\times m\\
\vdots & \vdots\\
{\huge 0}_{n\times\left[  \left(  n-1\right)  m\right]  } & \left.
\begin{array}
[c]{ccc}%
0 & \cdots & 0\\
\vdots & \ddots & \vdots\\
\lambda_{2} & \cdots & 0
\end{array}
\right\}  n\times m
\end{array}
\right)  _{mn\times mn},
\]
it is easy to see that for $l\leq-1$, $B_{0}^{\left(  l\right)  }$ is a
transition rate matrix from Level $l+1$ to Level $l$ while Level $l+1$ and
Level $l$ have the same lexicographic order, we obtain%
\[
B_{0}^{\left(  l\right)  }=\left(
\begin{array}
[c]{cc}%
\text{\ }{\huge 0}_{m\times\left[  \left(  n-1\right)  m\right]  }\text{\ } &
\left.
\begin{array}
[c]{ccc}%
\lambda_{2} &  & \\
& \ddots & \\
&  & \lambda_{2}%
\end{array}
\right\}  m\times m\\
{\huge 0}_{\left[  \left(  n-1\right)  m\right]  \times\left[  \left(
n-1\right)  m\right]  } & \text{ \ \ }{\huge 0}_{\left[  \left(  n-1\right)
m\right]  \times m}%
\end{array}
\right)  _{mn\times mn},\text{ }l\leq-1;
\]
since $B_{2}^{\left(  -1\right)  }$ is a transition rate matrix from Level
$-1$ to Level $0$ while Level $-1$ and Level $0$ have a different
lexicographic order, we obtain%
\[
B_{2}^{\left(  -1\right)  }=\left(
\begin{array}
[c]{cccc}%
\left.
\begin{array}
[c]{ccccc}%
0 & 0 & \cdots & 0 & \lambda_{1}\\
\vdots & \vdots &  & \vdots & \vdots\\
0 & 0 & \cdots & 0 & 0\\
0 & 0 & \cdots & 0 & 0
\end{array}
\right\}  m\times n &  &  & \\
\left.
\begin{array}
[c]{ccccc}%
0 & 0 & \cdots & \lambda_{1} & 0\\
\vdots & \vdots &  & \vdots & \vdots\\
0 & 0 & \cdots & 0 & 0\\
0 & 0 & \cdots & 0 & 0
\end{array}
\right\}  m\times n &  &  & \\
\vdots &  &  & \\
\left.
\begin{array}
[c]{ccccc}%
0 & \lambda_{1} & \cdots & 0 & 0\\
\vdots & \vdots &  & \vdots & \vdots\\
0 & 0 & \cdots & 0 & 0\\
0 & 0 & \cdots & 0 & 0
\end{array}
\right\}  m\times n &  &  & \\
\left.
\begin{array}
[c]{ccccc}%
\lambda_{1} & 0 & \cdots & 0 & 0\\
\vdots & \vdots &  & \vdots & \vdots\\
0 & 0 & \cdots & 0 & 0\\
0 & 0 & \cdots & 0 & n\theta_{2}%
\end{array}
\right\}  m\times n & \left.
\begin{array}
[c]{cc}%
{\huge 0} & 0\\
\vdots & \vdots\\
{\huge 0} & n\theta_{2}\\
{\huge 0} & 0
\end{array}
\right\}  m\times n & \cdots & \left.
\begin{array}
[c]{cc}%
{\huge 0} & n\theta_{2}\\
\vdots & \vdots\\
{\huge 0} & 0\\
{\huge 0} & 0
\end{array}
\right\}  m\times n
\end{array}
\right)  _{_{mn\times mn}},
\]
observing that for $l\leq-1$, $B_{2}^{\left(  l\right)  }$ is a transition
rate matrix from Level $l$ to Level $l+1$ while Level $l$ and Level $l+1$ have
the same lexicographic order, we obtain%
\[
B_{2}^{\left(  l\right)  }=\left(
\begin{array}
[c]{ccc}%
\left.
\begin{array}
[c]{ccc}%
0 &  & \lambda_{1}\\
& \ddots & \\
&  & 0
\end{array}
\right\}  m\times m &  & \\
\vdots & \ddots & \\
\left.
\begin{array}
[c]{ccc}%
-ln\theta_{2} &  & \\
& \ddots & \\
&  & -ln\theta_{2}%
\end{array}
\right\}  m\times m & \cdots & \left.
\begin{array}
[c]{ccc}%
0 &  & \lambda_{1}\\
& \ddots & \\
&  & 0
\end{array}
\right\}  m\times m
\end{array}
\right)  _{mn\times mn},\text{ }l\leq-2,
\]

\begin{Rem}
The element order of the state space $\Omega=\cup_{k=-\infty}^{\infty}$Level
$k$ is arranged from left to right not only for the levels but also for the
elements in each level. Thus, from (\ref{L0}), (\ref{Lk}) and (\ref{Ll}), it
is observed the useful difference between Level $k$ for $k\geq0$, and Level
$L$ for $l\leq-1$. In this case, from considering the example in Figure 1, it
is useful to observe the two special blocks for how to write the infinitesiaml
generator $Q$ as follows:%
\[
B_{0}^{\left(  0\right)  }=%
\begin{array}
[c]{cc}
&
\begin{array}
[c]{cccccc}%
\left(  1,5\right)  & \left(  0,5\right)  & \left(  1,4\right)  & \left(
0,4\right)  & \left(  1,3\right)  & \left(  0,3\right)
\end{array}
\\%
\begin{array}
[c]{c}%
\left(  0,0\right) \\
\left(  0,1\right) \\
\left(  0,2\right) \\
\left(  1,0\right) \\
\left(  1,1\right) \\
\left(  1,2\right)
\end{array}
& \text{Level }0\rightarrow\text{Level }-1
\end{array}
\]
and%
\[
B_{2}^{\left(  -1\right)  }=%
\begin{array}
[c]{cc}
&
\begin{array}
[c]{cccccc}%
\left(  0,0\right)  & \left(  0,1\right)  & \left(  0,2\right)  & \left(
1,0\right)  & \left(  1,1\right)  & \left(  1,2\right)
\end{array}
\\%
\begin{array}
[c]{c}%
\left(  1,5\right) \\
\left(  0,5\right) \\
\left(  1,4\right) \\
\left(  0,4\right) \\
\left(  1,3\right) \\
\left(  0,3\right)
\end{array}
& \text{Level }-1\rightarrow\text{Level }0
\end{array}
.
\]
\end{Rem}

Now, we discuss the stability of the QBD process $Q$ with bidirectional
infinite levels. Our method is to divide the bidirectional QBD process $Q$
into two unilateral QBD processes: $Q_{A}$ and $Q_{B}$.

To analyze the matched queue with matching batch pair $(m,n)$, it is
worthwhile to note that a simple relation between Levels $0$ and $-1$ is a
key. From Levels $0$ and $-1$, we can divide the bidirectional QBD process $Q$
into two unilateral QBD processes: $Q_{A}$ and $Q_{B}$. Based on this, the
infinitesimal generators of the two unilateral QBD processes $Q_{A}$ and
$Q_{B}$ are respectively given by
\[
Q_{A}=\left(
\begin{array}
[c]{cccccc}%
C+B_{0}^{\left(  0\right)  } & \fbox{$A_{0}^{\left(  0\right)  }$} &  &  &  &
\\
\fbox{$A_{2}^{\left(  1\right)  }$} & A_{1}^{\left(  1\right)  } &
A_{0}^{\left(  1\right)  } &  &  & \\
& A_{2}^{\left(  2\right)  } & A_{1}^{\left(  2\right)  } & A_{0}^{\left(
2\right)  } &  & \\
&  & A_{2}^{\left(  3\right)  } & A_{1}^{\left(  3\right)  } & A_{0}^{\left(
3\right)  } & \\
&  &  & \ddots & \ddots & \ddots
\end{array}
\right)
\]
and%
\[
Q_{B}=\left(
\begin{array}
[c]{cccccc}%
C+A_{0}^{\left(  0\right)  } & \fbox{$B_{0}^{\left(  0\right)  }$} &  &  &  &
\\
\fbox{$B_{2}^{\left(  -1\right)  }$} & B_{1}^{\left(  -1\right)  } &
B_{0}^{\left(  -1\right)  } &  &  & \\
& B_{2}^{\left(  -2\right)  } & B_{1}^{\left(  -2\right)  } & B_{0}^{\left(
-2\right)  } &  & \\
&  & B_{2}^{\left(  -3\right)  } & B_{1}^{\left(  -3\right)  } &
B_{0}^{\left(  -3\right)  } & \\
&  &  & \ddots & \ddots & \ddots
\end{array}
\right)  ,
\]
which can be given by the following infinitesimal generator when the levels
are re-arranged in a new order $\left\{  \text{Level }0,\text{Level
}-1,\text{Level }-2,\ldots\right\}  $ from its original order $\{\ldots,$Level
$-2$, Level $-1,$Level $0\}$, seeing the above part of the infinitesimal
generator $Q$ given in (\ref{Q}) as follows:%
\[
Q_{B}=\left(
\begin{array}
[c]{cccccc}%
\ddots & \ddots & \ddots &  &  & \\
& B_{0}^{\left(  -3\right)  } & B_{1}^{\left(  -3\right)  } & B_{2}^{\left(
-3\right)  } &  & \\
&  & B_{0}^{\left(  -2\right)  } & B_{1}^{\left(  -2\right)  } &
B_{2}^{\left(  -2\right)  } & \\
&  &  & B_{0}^{\left(  -1\right)  } & B_{1}^{\left(  -1\right)  } &
\fbox{$B_{2}^{\left(  -1\right)  }$}\\
&  &  &  & \fbox{$B_{0}^{\left(  0\right)  }$} & C+A_{0}^{\left(  0\right)  }%
\end{array}
\right)  .
\]
Note that $\left[  \left(  C+B_{0}^{\left(  0\right)  }\right)  +A_{0}%
^{\left(  0\right)  }\right]  e=0$\ for $Q_{A}$, and $\left[  \left(
C+A_{0}^{\left(  0\right)  }\right)  +B_{0}^{\left(  0\right)  }\right]
e=0$\ for $Q_{B}$, where $e$ is a column vector of ones with a suitable dimension.

\vskip                                                           0.4cm

In the remainder of this section, we study the stability of the matched queue
with matching batch pair $(m,n)$. It is easy to see that the customers'
impatient behavior plays a key role in the stability analysis.

The following theorem provides a sufficient condition under which the QBD
process with bidirectional infinite levels is stable.

\begin{The}
\ \label{The1}If $\left(  \theta_{1},\theta_{2}\right)  >0$, then the QBD
process $Q$ with bidirectional infinite levels is irreducible and positive
recurrent. Thus, the matched queue with matching batch pair $(m,n)$ is stable.
\end{The}

\noindent\textbf{Proof.}\textit{ }Please see Proof of Theorem \ref{The1} in
the appendix.

\begin{Rem}
(a) To develop the matrix-analytic method of matched queues, we have to assume
that the impatient times are exponential. Even so, it is still very
complicated to write the infinitesimal generator $Q$ given in (\ref{Q}).

(b) Similarly to the retrial times in a retrial queue, the non-exponential
impatient times can make a Markov modeling very complicated due to the
parallel work of multiple impatient times. If $m$ impatient times are of phase
type, then so far we have not known how to write the infinitesimal generator
$Q$ yet. Further, if $m$ impatient times have general distributions, it
becomes more difficult due to the fact that an embedded Markov chain has to be
established accordingly.
\end{Rem}

\section{The Stationary Queue Length}

In this section, we first provide a bilateral matrix-product expression for
the stationary probability vector of the level-dependent QBD process with
bidirectional infinite levels by means of the RG-factorizations given in Li
\cite{Li:2010}. Then we compute the two average stationary queue length for
any A- or B-customer.

We write%
\[
p_{i,j}\left(  t\right)  =P\left\{  N_{1}\left(  t\right)  =i\text{, }%
N_{2}\left(  t\right)  =j\right\}  .
\]
Since the level-dependent QBD process with bidirectional infinite levels is
stable, we have%
\[
\pi_{i,j}=\lim_{t\rightarrow+\infty}p_{i,j}\left(  t\right)  .
\]
For $l=-1,-2,-3,...$, we write%
\begin{align*}
\pi_{l}  &  =\left(  \pi_{m-1,\left(  -l\right)  n+n-1},\pi_{m-2,\left(
-l\right)  n+n-1},\ldots,\pi_{0,\left(  -l\right)  n+n-1};\ldots;\right. \\
&  \text{ \ \ \ }\left.  \pi_{m-1,\left(  -l\right)  n},\pi_{m-2,\left(
-l\right)  n},\ldots,\pi_{0,\left(  -l\right)  n}\right)  ,
\end{align*}
for $k=0$, we write%
\[
\pi_{0}=\left(  \pi_{0,0},\pi_{0,1},\ldots,\pi_{0,n-1};\pi_{1,0},\pi
_{1,1},\ldots,\pi_{1,n-1};\ldots;\pi_{m-1,0},\pi_{m-1,1},\ldots,\pi
_{m-1,n-1}\right)  ,
\]
for $k=1,2,3,...$, we write%
\begin{align*}
\pi_{k}  &  =\left(  \pi_{km,0},\pi_{km,1},\ldots,\pi_{km,n-1};\pi
_{km+1,0},\pi_{km+1,1},\ldots,\pi_{km+1,n-1};\ldots;\right. \\
&  \text{ \ \ \ }\left.  \pi_{km+m-1,0},\pi_{km+m-1,1},\ldots,\pi
_{km+m-1,n-1}\right)  ,
\end{align*}
and%
\[
\pi=\left(  \ldots,\pi_{-2},\pi_{-1},\pi_{0},\pi_{1},\pi_{2},\ldots\right)  .
\]

To compute the stationary probability vector of the level-dependent QBD
process $Q$ with bidirectional infinite levels, we first need to compute the
stationary probability vectors of the two unilateral QBD processes $Q_{A}$ and
$Q_{B}$. Then we use Levels $1$, $0$ and $-1$ as three interaction boundary
levels, which are used to further determine the stationary probability vectors.

Note that the two unilateral QBD processes $Q_{A}$ and $Q_{B}$ are
level-dependent, thus we need to apply the RG-factorization given in Li
\cite{Li:2010} to calculate their stationary probability vectors. To this end,
we need to introduce the $U$-, $R$\textit{-} and $G$\textit{-}measures for the
two unilateral QBD processes $Q_{A}$ and $Q_{B}$, respectively. In fact, such
a level-dependent QBD process was analyzed in Li and Cao \cite{Li:2004}.

For the unilateral QBD process $Q_{A}$, we define the $UL$-type $U$-,
$R$\textit{-} and $G$\textit{-}measures as%
\begin{align}
U_{0}  &  =\left(  C+B_{0}^{\left(  0\right)  }\right)  +A_{0}^{\left(
0\right)  }\left(  -U_{1}^{-1}\right)  A_{2}^{\left(  1\right)  },\text{
\ }\nonumber\\
U_{k}  &  =A_{1}^{\left(  k\right)  }+A_{0}^{\left(  k\right)  }\left(
-U_{k+1}^{-1}\right)  A_{2}^{\left(  k+1\right)  },\text{ \ }k\geq1
\label{equ2}%
\end{align}%
\[
R_{k}=A_{0}^{\left(  k\right)  }\left(  -U_{k+1}^{-1}\right)  ,\text{ \ }%
k\geq0,
\]
and%
\[
G_{k}=\left(  -U_{k}^{-1}\right)  A_{2}^{\left(  k\right)  },\text{ \ }%
k\geq1.
\]
Obviously, it is well-known from Li and Cao \cite{Li:2004} that the matrix
sequence $\left\{  R_{k},k\geq0\right\}  $ is the minimal nonnegative solution
to the system of nonlinear matrix equations%

\begin{equation}
A_{0}^{\left(  k\right)  }+R_{k}A_{1}^{\left(  k+1\right)  }+R_{k}R_{k+1}%
A_{2}^{\left(  k+2\right)  }=0,\text{ \ }k\geq0; \label{equ3}%
\end{equation}
and the matrix sequence $\left\{  G_{k},k\geq1\right\}  $ is the minimal
nonnegative solution to the system of nonlinear matrix equations%
\begin{equation}
A_{0}^{\left(  k\right)  }G_{k+1}G_{k}+A_{1}^{\left(  k\right)  }G_{k}%
+A_{2}^{\left(  k\right)  }=0,\text{ \ }k\geq1. \label{equ4}%
\end{equation}
It is worthwhile to note that the systems (\ref{equ3}) and (\ref{equ4}) of
nonlinear matrix equations were first given in Ramaswami and Taylor
\cite{Ram:1996}.

Let $\pi_{A}=\left(  \pi_{0}^{A},\pi_{1}^{A},\pi_{2}^{A},\ldots\right)  $ be
the stationary probability vector of the unilateral QBD process $Q_{A}$. Then
from Subsection 2.7.3 in Chapter 2 of Li \cite{Li:2010} or Li and Cao Li
\cite{Li:2004}, by using the $R$-measure $\left\{  R_{k}:k\geq0\right\}  $ we
have
\begin{equation}
\pi_{k}^{A}=\gamma^{A}\pi_{0}^{A}R_{0}R_{1}\cdots R_{k-1},\text{ \ }k\geq1,
\label{equ6}%
\end{equation}
where $\pi_{0}^{A}$ is the stationary probability vector of the censored chain
$U_{0}$ to level $0$, and $\gamma^{A}$ is a regularization coefficient such
that $\sum_{k=0}^{\infty}\pi_{k}^{A}e=1$. Note that the expression
(\ref{equ6}) of the stationary probability vector was first obtained by
Ramaswami and Taylor \cite{Ram:1996}.

By conducting a similar analysis to those in Equations (\ref{equ2}) to
(\ref{equ6}), we can give the stationary probability vector, $\pi_{B}=\left(
\pi_{0}^{B},\pi_{1}^{B},\pi_{2}^{B},\ldots\right)  $, of the unilateral QBD
process $Q_{B}$. Here, we only provide the $R$-measure $\left\{
\mathbb{R}_{l}:l\leq0\right\}  $, while the $U$-measure $\left\{
\mathbb{U}_{l}:l\leq0\right\}  $ and $G$-measure $\left\{  \mathbb{G}%
_{l}:l\leq-1\right\}  $ is omitted for brevity.

Let the matrix sequence $\left\{  \mathbb{R}_{l},l\leq0\right\}  $ be the
minimal nonnegative solution to the system of nonlinear matrix equations%

\begin{equation}
B_{0}^{\left(  l\right)  }+\mathbb{R}_{l}B_{1}^{\left(  l-1\right)
}+\mathbb{R}_{l}\mathbb{R}_{l-1}B_{2}^{\left(  l-2\right)  }=0,\text{ }l\leq0,
\label{equ7}%
\end{equation}
By using the $R$-measure $\left\{  \mathbb{R}_{l}:l\leq0\right\}  $ we obtain
\begin{equation}
\pi_{l}^{B}=\gamma^{B}\pi_{0}^{B}\mathbb{R}_{0}\mathbb{R}_{-1}\cdots
\mathbb{R}_{l+1},\text{ \ }l\leq-1, \label{equ8}%
\end{equation}
where $\pi_{0}^{B}$ is the stationary probability vector of the censored chain
$\mathbb{U}_{0}$ to level $0$, and $\gamma^{B}$ is a regularization
coefficient such that $\sum_{k=0}^{\infty}\pi_{k}^{B}e=1$.

The following theorem expresses the stationary probability vector $\pi=\left(
\ldots,\pi_{-2},\pi_{-1},\right.  $

\noindent$\left.  \pi_{0},\pi_{1},\pi_{2},\ldots\right)  $ of the
level-dependent QBD process $Q$ with bidirectional infinite levels by means of
the stationary probability vectors $\pi_{A}=\left(  \pi_{2}^{A},\pi_{3}%
^{A},\pi_{4}^{A},\ldots\right)  $ given in (\ref{equ6}), and $\pi_{B}=\left(
\pi_{-2}^{B},\pi_{-3}^{B},\pi_{-4}^{B},\ldots\right)  $ given in (\ref{equ8}).

\begin{The}
\label{The2}The stationary probability vector $\pi$ of the level-dependent QBD
process $Q$ with bidirectional infinite levels is given by
\begin{equation}
\pi_{k}=c\widetilde{\pi}_{k}, \label{equ10}%
\end{equation}
and%
\begin{equation}
\widetilde{\pi}_{k}=\left\{
\begin{array}
[c]{l}%
\widetilde{\pi}_{-1}\mathbb{R}_{-1}\mathbb{R}_{-2}\cdots\mathbb{R}%
_{k+1},\text{ }k\leq-2,\\
\widetilde{\pi}_{-1},\\
\widetilde{\pi}_{0},\\
\widetilde{\pi}_{1},\\
\widetilde{\pi}_{1}R_{1}R_{2}\cdots R_{k-1},\text{ \ }k\geq2,
\end{array}
\right.  \label{equ11}%
\end{equation}
where the three boundary vectors $\widetilde{\pi}_{-1},\widetilde{\pi}%
_{0},\widetilde{\pi}_{1}$ are uniquely determined by the following system of
linear equations%
\begin{equation}
\left\{
\begin{array}
[c]{l}%
\widetilde{\pi}_{0}A_{0}^{\left(  0\right)  }+\widetilde{\pi}_{1}\left[
A_{1}^{\left(  1\right)  }+R_{1}A_{2}^{\left(  2\right)  }\right]  =0,\\
\pi_{-1}B_{2}^{\left(  -1\right)  }+\widetilde{\pi}_{0}C+\widetilde{\pi}%
_{1}A_{2}^{\left(  1\right)  }=0,\\
\widetilde{\pi}_{0}B_{0}^{\left(  0\right)  }+\widetilde{\pi}_{-1}\left[
B_{1}^{\left(  -1\right)  }+\mathbb{R}_{-1}B_{2}^{\left(  -2\right)  }\right]
=0,
\end{array}
\right.  \label{equ12}%
\end{equation}
and the positive constant $c$ is uniquely given by%
\begin{equation}
c=\frac{1}{\sum_{k\leq-2}\widetilde{\pi}_{-1}\mathbb{R}_{-1}\mathbb{R}%
_{-2}\cdots\mathbb{R}_{k+1}\mathbf{e}+\widetilde{\pi}_{-1}\mathbf{e}%
+\widetilde{\pi}_{0}\mathbf{e}+\widetilde{\pi}_{1}\mathbf{e}+\sum
_{k=2}^{\infty}\widetilde{\pi}_{1}R_{1}R_{2}\cdots R_{k-1}\mathbf{e}}.
\label{equ121}%
\end{equation}
\end{The}

\noindent\textbf{Proof.}\textit{ }Please see Proof of Theorem \ref{The2} in
the appendix.

\begin{Rem}
As seen from Theorem \ref{The2}, there is no explicit expression for the
stationary probability vector of the level-dependent QBD process with
bidirectional infinite levels, in which the impatient customers lead to a
level-dependent Markov process whose stationary probability vector computation
is more complicated by means of the RG-factorizations. Thus there does not
exist an analytic expression which is use to further discuss performance
measures of the matched queue with matching batch pair $(m,n)$.
\end{Rem}

In the remainder of this section, we compute the two average stationary queue
lengths for the A- and B-customers, respectively.

Note that the matched queue with matching batch pair $(m,n)$ is stable for
$\left(  \theta_{1},\theta_{2}\right)  >0$, we denote by $\mathcal{Q}^{\left(
1\right)  }$ and $\mathcal{Q}^{\left(  2\right)  }$ the stationary queue
lengths of the A- and B-customers, respectively. By using Theorem \ref{The2},
we provide the average stationary queue lengths of the A- and B-customers as follows:

\textbf{(a) }The average stationary queue length of the A-customers is given by%

\[
E\left[  \mathcal{Q}^{\left(  1\right)  }\right]  =\sum_{k=1}^{\infty}%
\sum_{i=0}^{m-1}\left(  km+i\right)  \sum_{j=0}^{n-1}\pi_{km+i,j}%
+\sum_{-\infty<l\leq0}\sum_{i=0}^{m-1}i\sum_{j=0}^{n-1}\pi_{i,\left(
-l\right)  n+j}.
\]

\textbf{(b) }The average stationary queue length of the B-customers is given by%

\[
E\left[  \mathcal{Q}^{\left(  2\right)  }\right]  =\sum_{k=0}^{\infty}%
\sum_{j=0}^{n-1}j\sum_{i=0}^{m-1}\pi_{km+i,j}+\sum_{-\infty<l\leq-1}\sum
_{j=0}^{n-1}\left[  \left(  -l\right)  n+j\right]  \sum_{i=0}^{m-1}%
\pi_{i,\left(  -l\right)  n+j}.
\]

By using the efficient algorithms given in Bright and Taylor \cite{Bri:1995,
Bri:1997} (e.g, see Liu et al. \cite{Liu:2020} for more details), we conduct
numerical examples to analyze how the average stationary queue lengths of A-
and B-customers are influenced by two key parameters: $\theta_{1}$ and
$\theta_{2}$. To this end, we take the system parameters: $\lambda_{1}=1$,
$\lambda_{2}=2$, $m=2$ and $n=3$ for the purpose of illustration.

Figure 3 shows that the average stationary queue length $E\left[
\mathcal{Q}^{\left(  1\right)  }\right]  $ decreases as $\theta_{1}$
increases, while it increases as $\theta_{2}$ increases.

From Figure 4, it is easy seen that the average stationary queue length
$E\left[  \mathcal{Q}^{\left(  2\right)  }\right]  $ increases as $\theta_{1}$
increases, while it decreases as $\theta_{2}$ increases.

\begin{figure}[h]
\centering  \includegraphics[width=7cm]{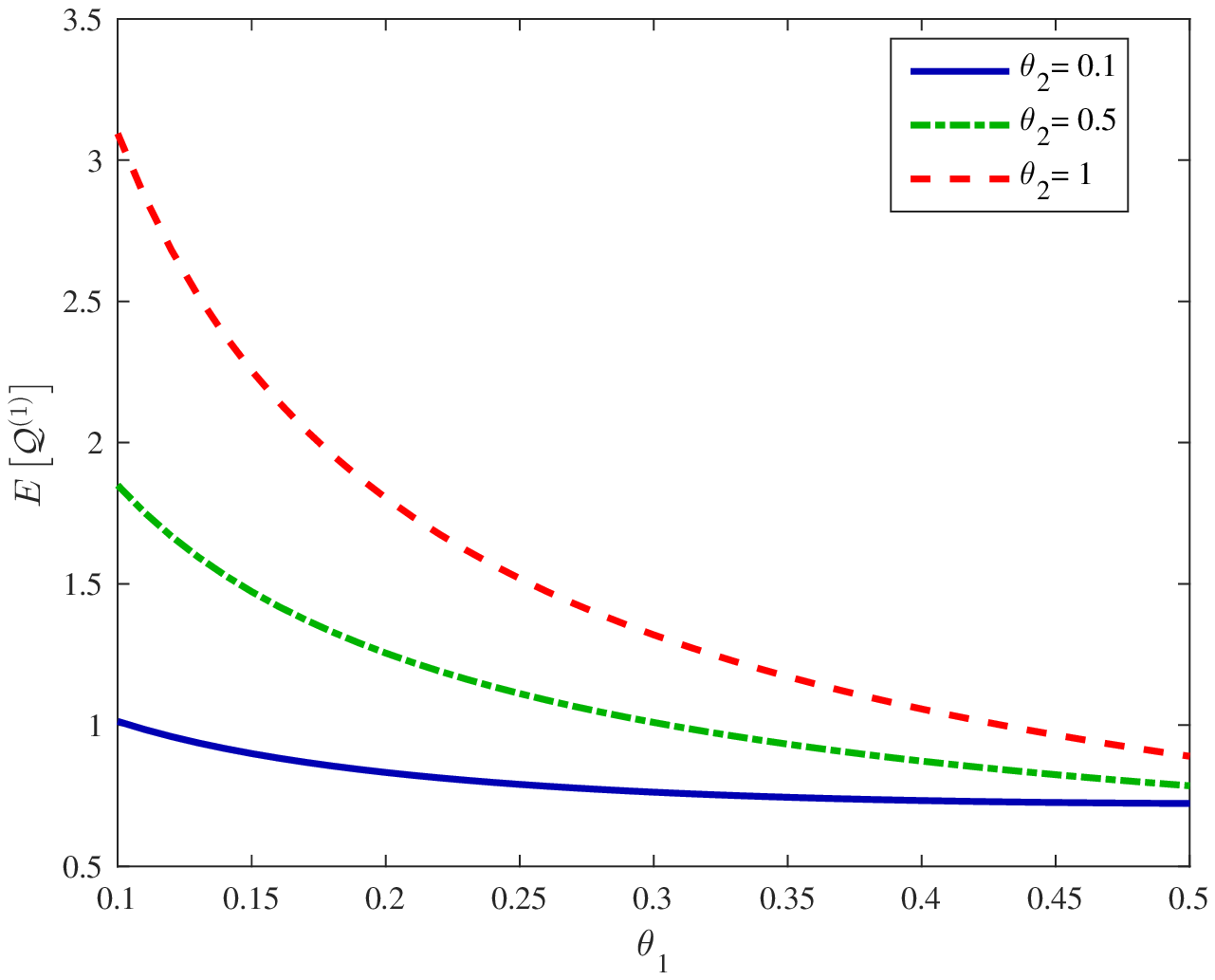}  \centering
\includegraphics[width=7cm]{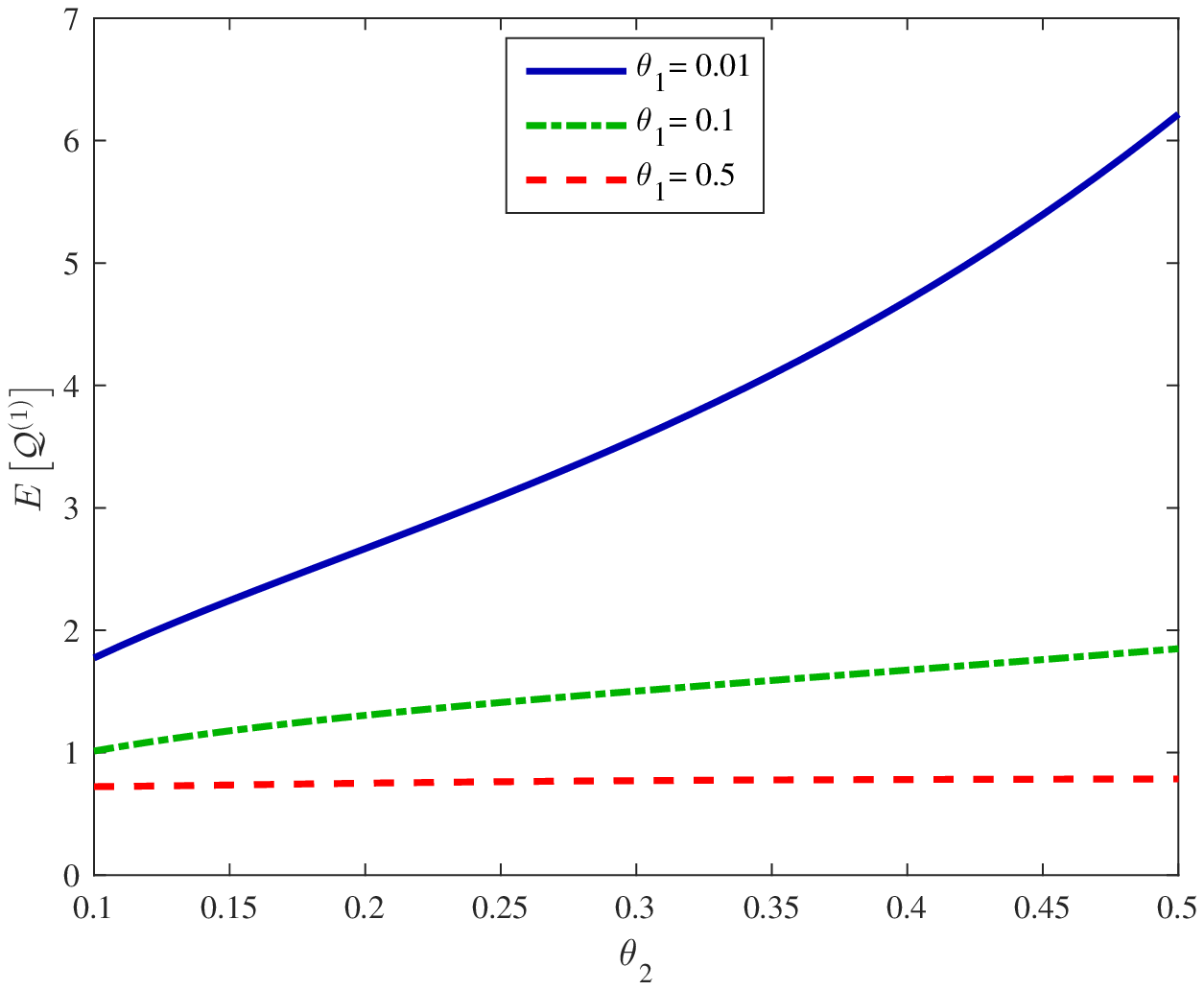}  \caption{$E\left[  \mathcal{Q}%
^{\left(  1\right)  }\right]  $ vs $\theta_{1}$ and $\theta_{2}$}%
\label{Fig-3}%
\end{figure}

\begin{figure}[h]
\centering  \includegraphics[width=7cm]{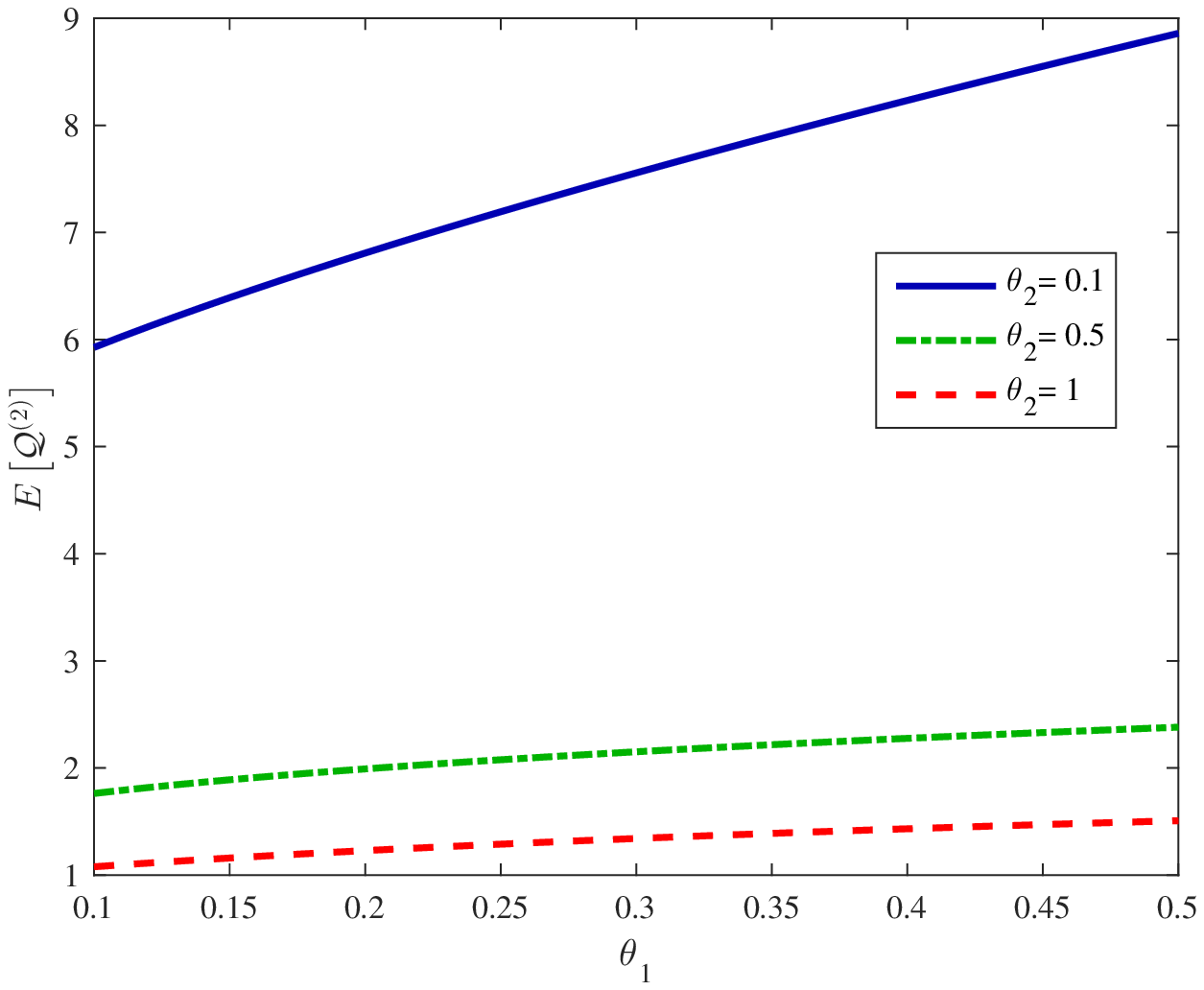}  \centering
\includegraphics[width=7cm]{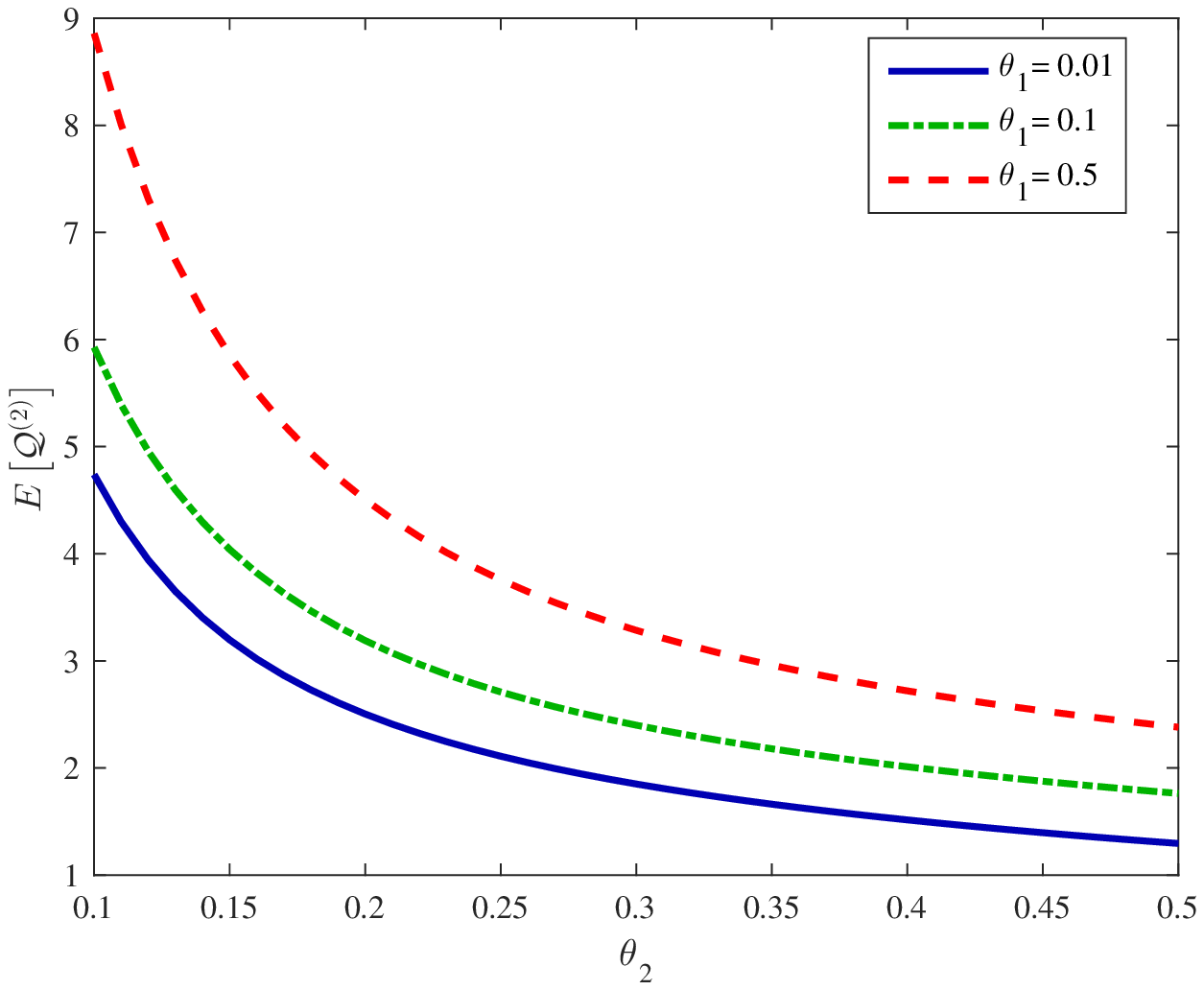}  \caption{$E\left[  \mathcal{Q}%
^{\left(  2\right)  }\right]  $ vs $\theta_{1}$ and $\theta_{2}$}%
\label{Fig-4}%
\end{figure}

The two numerical results are intuitive. As $\theta_{1}$ increases, more and
more A-customers quickly leave the system so that $E\left[  \mathcal{Q}%
^{\left(  1\right)  }\right]  $ decreases. On the other hand, as $\theta_{2}$
increases, more and more B-customers quickly leave the system so that the
probability that an A-customer can match a B-customer will become smaller and
smaller. Thus, $E\left[  \mathcal{Q}^{\left(  1\right)  }\right]  $ increases
as $\theta_{2}$ increases.

\section{The Sojourn Time}

In this section, we compute the average sojourn time of any arriving A- or
B-customer in the matched queue with matching batch pair $(m,n)$. Our analysis
includes three different parts: Using the Little's formula, a probabilistic
calculation, and an upper bound. Based on this, we can further find some
useful random relations in this matched queue.

\textbf{Part one: Using the Little's formula}

When the matched queue is stable, we denote by $W$ the sojourn time of any
arriving A-customer. By using the Little's formula and the average stationary
queue length $E\left[  \mathcal{Q}^{\left(  1\right)  }\right]  $ of the
A-customers, we obtain%
\begin{equation}
E\left[  W\right]  =\frac{E\left[  \mathcal{Q}^{\left(  1\right)  }\right]
}{\lambda_{1}}. \label{S-0}%
\end{equation}

\textbf{Part two: A probabilistic calculation}

Although the Little's formula provides a simple method to compute the average
sojourn time $E\left[  W\right]  $, it is still necessary and useful for
discussing the random structure of $E\left[  W\right]  $. To this end, our
analysis contains three different cases as follows:

(a) An arriving A-customer observes the system state $\left(  i,j\right)  $
for $0\leq i\leq m-1$ and $j\geq n$;

(b) an arriving A-customer observes the system state $\left(  i,j\right)  $
for $0\leq i\leq m-1$ and $0\leq j\leq n-1$; and

(c) an arriving A-customer observes the system state $\left(  i,j\right)  $
for $i\geq m$ and $0\leq j\leq n-1$.

The following lemma is useful for our computation in the above three cases,
while its proof is easy and is omitted here.

\begin{Lem}
\label{Lem:Soj}Let $C$ and $D$ be two independent nonnegative random
variables. Then%
\[
E\left[  \min\left\{  C,D\right\}  \right]  =\int_{0}^{+\infty}P\left\{
\min\left\{  C,D\right\}  >x\right\}  \text{d}x=\int_{0}^{+\infty}P\left\{
C>x\right\}  P\left\{  C>x\right\}  \text{d}x.
\]
\end{Lem}

Now, we use Lemma \ref{Lem:Soj} to compute the average sojourn time $E\left[
W\right]  $ from three different cases.

\textbf{Case a: An arriving A-customer observes the system state }$\left(
i,j\right)  $\textbf{ for }$0\leq i\leq m-1$\textbf{ and }$j\geq n$

In this case, if $i=m-1$, then the arriving A-customer makes the number of
A-customers become $m$. In this case, the $m$ A-customers can match $n$
B-customers as a group, the $m+n$ customers leave this system immediately.
Thus the average sojourn time is given by%
\begin{equation}
E\left[  W\right]  =0. \label{S-1}%
\end{equation}

If $0\leq i\leq m-2$, then the arriving A-customer still needs to wait for
$m-\left(  i+1\right)  $ arrivals of A-customers such that the number of
A-customers is $m$. Thus the $m$ A-customers can match $n$ B-customers as a
group, which leaves this system immediately.

Let $X_{k}$ be the $k$th interarrival time of the Poisson process with arrival
rate $\lambda_{1}$, and $Z$ be the exponential impatient time of A-customer
with impatient rate $\theta_{1}$. Then%
\[
W=\min\left\{  Z,\sum_{k=1}^{m-\left(  i+1\right)  }X_{k}\right\}  .
\]
This gives%
\[
E\left[  W\right]  =\int_{0}^{+\infty}P\left\{  Z>x\right\}  P\left\{
\sum_{k=1}^{m-\left(  i+1\right)  }X_{k}>x\right\}  \text{d}x.
\]
Note that $\sum_{k=1}^{m-\left(  i+1\right)  }X_{k}$ is an Erlang-$\left(
m-i-1\right)  $ distribution, we have%
\[
P\left\{  \sum_{k=1}^{m-\left(  i+1\right)  }X_{k}>u\right\}  =e^{-\lambda
_{1}u}\sum_{k=0}^{m-\left(  i+1\right)  -1}\frac{\left(  \lambda_{1}u\right)
^{k}}{k!},
\]
this gives%
\begin{equation}
E\left[  W\right]  =\int_{0}^{+\infty}e^{-\left(  \lambda_{1}+\theta
_{1}\right)  x}\sum_{k=0}^{m-\left(  i+1\right)  -1}\frac{\left(  \lambda
_{1}x\right)  ^{k}}{k!}\text{d}x. \label{S-2}%
\end{equation}

\textbf{Case b: An arriving A-customer observes the system state }$\left(
i,j\right)  $\textbf{ for }$0\leq i\leq m-1$ \textbf{and} $0\leq j\leq n-1$

In this case, if $i=m-1$, then the arriving A-customer makes that the number
of A-customers become $m$. However, they still need to wait for $n-j$ arrivals
of B-customers such that the number of B-customers is $n$. In this case, the
$n$ B-customers can match $m$ A-customers as a group, which leaves this system immediately.

Let $Y_{k}$ be the $k$th interarrival time of the Poisson process with arrival
rate $\lambda_{2}$. Then%
\[
W=\min\left\{  Z,\sum_{k=1}^{n-j}Y_{k}\right\}  .
\]
This gives%
\begin{align}
E\left[  W\right]   &  =\int_{0}^{+\infty}P\left\{  Z>x\right\}  P\left\{
\sum_{k=1}^{n-j}Y_{k}>x\right\}  \text{d}x\nonumber\\
&  =\int_{0}^{+\infty}e^{-\left(  \lambda_{2}+\theta_{1}\right)  x}\sum
_{k=0}^{n-j-1}\frac{\left(  \lambda_{2}x\right)  ^{k}}{k!}\text{d}x.
\label{S-3}%
\end{align}

If $0\leq i\leq m-2$, then we still need to not only wait for $m-\left(
i+1\right)  $ arrivals of A-customers such that the number of A-customers is
$m$, but also wait for $n-j$ arrivals of B-customers such that the number of
B-customers is $n$. In this case, the $n$ B-customers can match $m$
A-customers as a group, which leaves this system immediately. Thus we obtain%
\[
W=\min\left\{  Z,\max\left\{  \sum_{k=1}^{m-\left(  i+1\right)  }X_{k}%
,\sum_{k=1}^{n-j}Y_{k}\right\}  \right\}  .
\]
This gives%
\begin{align}
E\left[  W\right]   &  =\int_{0}^{+\infty}P\left\{  Z>x\right\}  P\left\{
\max\left\{  \sum_{k=1}^{m-\left(  i+1\right)  }X_{k},\sum_{k=1}^{n-j}%
Y_{k}\right\}  >x\right\}  \text{d}x\nonumber\\
&  =\int_{0}^{+\infty}P\left\{  Z>x\right\}  \left[  1-P\left\{  \max\left\{
\sum_{k=1}^{m-\left(  i+1\right)  }X_{k},\sum_{k=1}^{n-j}Y_{k}\right\}  \leq
x\right\}  \right]  \text{d}x\nonumber\\
&  =\int_{0}^{+\infty}P\left\{  Z>x\right\}  \left[  1-P\left\{  \sum
_{k=1}^{m-\left(  i+1\right)  }X_{k}\leq x\right\}  P\left\{  \sum_{k=1}%
^{n-j}Y_{k}\leq x\right\}  \right]  \text{d}x. \label{S-4}%
\end{align}

\textbf{Case c: An arriving A-customer observes the system state }$\left(
i,j\right)  $\textbf{ for }$i\geq m$ \textbf{and} $0\leq j\leq n-1$.

Since $i\geq m$, there exists a unique positive integer $h$ such that
$i=hm+f$, where $0\leq f\leq m-1$.

If $f=m-1$, then the arriving A-customer makes the number of A-customers
become $\left(  h+1\right)  m$. We still need to wait for $hn+n-j$ arrivals of
B-customers such that the number of B-customers becomes $\left(  h+1\right)
n$. In this case, the $\left(  h+1\right)  n$ B-customers can match $\left(
h+1\right)  m$ A-customers as $h+1$ groups, $\left(  h+1\right)  \left(
m+n\right)  $ customers leave this system immediately. Thus we obtain%
\[
W=\min\left\{  Z,\sum_{k=1}^{hn+n-j}Y_{k}\right\}  .
\]
This gives%
\begin{align}
E\left[  W\right]   &  =\int_{0}^{+\infty}P\left\{  Z>x\right\}  P\left\{
\sum_{k=1}^{hn+n-j}Y_{k}>x\right\}  \text{d}x\nonumber\\
&  =\int_{0}^{+\infty}e^{-\left(  \lambda_{2}+\theta_{1}\right)  x}\sum
_{k=0}^{hn+n-j-1}\frac{\left(  \lambda_{2}x\right)  ^{k}}{k!}\text{d}x.
\label{S-5}%
\end{align}

If $0\leq f\leq m-2$, then we still need to not only wait for $m-\left(
f+1\right)  $ arrivals of A-customers such that the number of A-customers
becomes $\left(  h+1\right)  m$, but also wait for $hn+n-j$ arrivals of
B-customers such that the number of B-customers becomes $\left(  h+1\right)
n$. In this case, the $\left(  h+1\right)  n$ B-customers can match $\left(
h+1\right)  m$ A-customers to form $h+1$ groups, which leave this system
immediately. Thus we obtain%
\[
W=\min\left\{  Z,\max\left\{  \sum_{k=1}^{m-f-1}X_{k},\sum_{k=1}^{hn+n-j}%
Y_{k}\right\}  \right\}  .
\]
This gives%
\begin{align}
E\left[  W\right]   &  =\int_{0}^{+\infty}P\left\{  Z>x\right\}  P\left\{
\max\left\{  \sum_{k=1}^{m-f-1}X_{k},\sum_{k=1}^{hn+n-j}Y_{k}\right\}
>x\right\}  \text{d}x\nonumber\\
&  =\int_{0}^{+\infty}P\left\{  Z>x\right\}  \left[  1-P\left\{  \sum
_{k=1}^{m-f-1}X_{k}\leq x\right\}  P\left\{  \sum_{k=1}^{hn+n-j}Y_{k}\leq
x\right\}  \right]  \text{d}x. \label{S-6}%
\end{align}

In what follows we provide an average stationary sojourn time of any arriving
A-customer by means of Equations (\ref{S-1}) to (\ref{S-6}) as well as the
stationary probability vector of this system. We write that for
$l=-1,-2,-3,...$,%
\begin{align*}
\Phi_{l}  &  =\left(  \phi_{m-1,\left(  -l\right)  n+n-1},\phi_{m-2,\left(
-l\right)  n+n-1},\ldots,\phi_{0,\left(  -l\right)  n+n-1};\ldots;\right. \\
&  \text{ \ \ \ }\left.  \phi_{m-1,\left(  -l\right)  n},\phi_{m-2,\left(
-l\right)  n},\ldots,\phi_{0,\left(  -l\right)  n}\right)  ,
\end{align*}
for $k=0$,%
\[
\Phi_{0}=\left(  \phi_{1,0},\phi_{1,1},\ldots,\phi_{1,n-1};\phi_{2,0}%
,\phi_{2,1},\ldots,\phi_{2,n-1};\ldots;\phi_{m-1,0},\phi_{m-1,1},\ldots
,\phi_{m-1,n-1}\right)  ,
\]
for $k=1,2,3,...$,%
\begin{align*}
\Phi_{k}  &  =\left(  \phi_{km,0},\phi_{km,1},\ldots,\phi_{km,n-1}%
;\phi_{km+1,0},\phi_{km+1,1},\ldots,\phi_{km+1,n-1};\ldots;\right. \\
&  \text{ \ \ \ }\left.  \phi_{km+m-1,0},\phi_{km+m-1,1},\ldots,\phi
_{km+m-1,n-1}\right)  ,
\end{align*}
and%
\begin{equation}
\Phi=\left(  \ldots,\Phi_{-2},\Phi_{-1},\Phi_{0},\Phi\pi_{1},\Phi_{2}%
,\ldots\right)  . \label{S-13}%
\end{equation}
For $1\leq i\leq m-1$ and $0\leq j\leq n-1$; $1\leq i\leq m-1$\textbf{
}and\textbf{ }$j\geq n$; and $i\geq m$ and $0\leq j\leq n-1$, we write
\[
\phi_{i,j}=\frac{1}{1-\sum\limits_{-\infty<l\leq0}\sum\limits_{j=0}^{n-1}%
\phi_{0,\left(  -l\right)  n+j}}\pi_{i,j}.
\]
It is clear that $\Phi e=1$.

When the matched queue with matching batch pair $(m,n)$ is stable, we assume
that an arriving A-customer enters State $\left(  i-1,j\right)  $ with
probability $\phi_{i,j}$ for $i\geq1$ at time $0$. Once the A-customer enters
the matched queue, the number of A-customers becomes $i$. In this case, we
obtain%
\begin{align}
E\left[  W\right]  =  &  \sum_{i=1}^{m-1}\sum_{j=n}^{\infty}\phi_{i,j}E\left[
W\text{ }|\text{ }\left(  N_{1}\left(  0\right)  ,N_{2}\left(  0\right)
\right)  =\left(  i,j\right)  \right]  \text{ \ \ \ \ \ \ \ \ \ \ \ Case
a}\nonumber\\
&  +\sum_{i=1}^{m-1}\sum_{j=0}^{n-1}\phi_{i,j}E\left[  W\text{ }|\text{
}\left(  N_{1}\left(  0\right)  ,N_{2}\left(  0\right)  \right)  =\left(
i,j\right)  \right]  \text{ \ \ \ \ \ \ \ \ Case b}\nonumber\\
&  +\sum_{i=m}^{\infty}\sum_{j=1}^{n-1}\phi_{i,j}E\left[  W\text{ }|\text{
}\left(  N_{1}\left(  0\right)  ,N_{2}\left(  0\right)  \right)  =\left(
i,j\right)  \right]  ,\text{\ \ \ \ \ \ \ \ Case c} \label{S-7}%
\end{align}
where $E\left[  W\text{ }|\text{ }\left(  N_{1}\left(  0\right)  ,N_{2}\left(
0\right)  \right)  =\left(  i,j\right)  \right]  $ is given in Case a by
(\ref{S-1}) and (\ref{S-2}); Case b by (\ref{S-3}) and (\ref{S-4}); and Case c
by (\ref{S-5}) and (\ref{S-6}).

\textbf{Part three: An upper bound}

Now, we provide a better upper bound of the average sojourn time $E\left[
W\right]  $ given in (\ref{S-7}) by means of a new PH distribution of
bidirectional infinite sizes.

We write a first passage time of the Markov process $Q$ (or $\left\{  \left(
N_{1}\left(  t\right)  ,N_{2}\left(  t\right)  \right)  ,\text{ }%
t\geq0\right\}  $) as%
\[
\xi=\inf\left\{  t:N_{1}\left(  t\right)  =0\right\}  ,
\]
that is, $\xi$ is such a first passage time that the the waiting room of
A-customers becomes empty for the first time. It is easy to see that $E\left[
W\right]  \leq E\left[  \xi\right]  $, since it is possible that the waiting
room of A-customers still contains some customers at the time that the
arriving A-customer leaves the system. Note that the arriving A-customer
enters State $\left(  i-1,j\right)  $ with probability $\phi_{i,j}$ for
$i\geq1$ at time $0$. It is also possible that the waiting room of A-customers
is empty at the time that the arriving A-customer leaves the system. Thus
$E\left[  \xi\right]  $ can be a better upper bound of the average sojourn
time $E\left[  W\right]  $.

Now, we compute the average first passage time $E\left[  \xi\right]  $. To do
this, we take all the states: $\left(  0,k\right)  $ for $k\geq0$, as an
absorbing state $\Delta$. Therefore, from the Markov process $Q$, we can set
up an absorbing Markov process whose infinitesmall generator is given by%
\[
\mathbb{Q}=\left(
\begin{array}
[c]{cc}%
0 & \mathbf{0}\\
T^{0} & T
\end{array}
\right)  ,
\]
where the first row and the first column of the matrix $\mathbb{Q}$ are
related to the absorbing state $\Delta$. To write the matrix $T$, we need to
check the levels of the absorbing Markov process $\mathbb{Q}$ as follows: For
$k\geq1,$
\[
\widetilde{\text{Level }k}=\text{Level }k,
\]%
\begin{align*}
\widetilde{\text{Level }0}=  &  \left\{  \left(  1,0\right)  ,\text{ }\left(
1,1\right)  ,,\ldots,\text{\ }\left(  1,n-1\right)  ;\left(  2,0\right)
,\left(  2,1\right)  ,\ldots,\left(  2,n-1\right)  ;\right. \\
&  \left.  \ldots;\text{\ }\left(  m-1,0\right)  ,\text{\ }\left(
m-1,1\right)  ,\ldots,\text{\ }\left(  m-1,n-1\right)  \right\}  ,
\end{align*}
and $l\leq-1,$%
\begin{align*}
\widetilde{\text{Level }l}=  &  \left\{  \left(  m-1,\left(  -l\right)
n+n-1\right)  ,\left(  m-2,\left(  -l\right)  n+n-1\right)  ,\ldots,\left(
1,\left(  -l\right)  n+n-1\right)  ;\ldots;\right. \\
&  \left(  m-1,\left(  -l\right)  n+1\right)  ,\left(  m-2,\left(  -l\right)
n+1\right)  ,\ldots,\left(  1,\left(  -l\right)  n+1\right)  ;\\
&  \left.  \left(  m-1,\left(  -l\right)  n\right)  ,\left(  m-2,\left(
-l\right)  n\right)  ,\ldots,\left(  1,\left(  -l\right)  n\right)  \right\}
,
\end{align*}
Thus the state space of the absorbing Markov process $\mathbb{Q}$ is given by%
\begin{align*}
\widetilde{\Omega}  &  =\left\{  \Delta\right\}  \cup\left\{  \bigcup
\limits_{-\infty<k<\infty}\widetilde{\text{Level }k}\right\} \\
&  =\left\{  \Delta\right\}  \cup\left\{  \bigcup\limits_{k=1}^{\infty
}\text{Level }k\right\}  \cup\left\{  \bigcup\limits_{-\infty<l\leq
0}\widetilde{\text{Level }l}\right\}  .
\end{align*}
By using the levels: Level $k$ for $k\geq1$ and $\widetilde{\text{Level }l}$
for $-\infty<l\leq0$, we obtain%
\[
T=\left(
\begin{array}
[c]{ccccccccccc}%
\ddots & \ddots & \ddots &  &  &  &  &  &  &  & \\
& \widetilde{B}_{0}^{\left(  -3\right)  } & \widetilde{B}_{1}^{\left(
-3\right)  } & \widetilde{B}_{2}^{\left(  -3\right)  } &  &  &  &  &  &  & \\
&  & \widetilde{B}_{0}^{\left(  -2\right)  } & \widetilde{B}_{1}^{\left(
-2\right)  } & \widetilde{B}_{2}^{\left(  -2\right)  } &  &  &  &  &  & \\
&  &  & \widetilde{B}_{0}^{\left(  -1\right)  } & \widetilde{B}_{1}^{\left(
-1\right)  } & \fbox{$\widetilde{B}_{2}^{\left(  -1\right)  }$} &  &  &  &  &
\\
&  &  &  & \fbox{$\widetilde{B}_{0}^{\left(  0\right)  }$} & \fbox{$\widetilde
{C}$} & \fbox{$\widetilde{A}_{0}^{\left(  0\right)  }$} &  &  &  & \\
&  &  &  &  & \fbox{$\widetilde{A}_{2}^{\left(  1\right)  }$} & A_{1}^{\left(
1\right)  } & A_{0}^{\left(  1\right)  } &  &  & \\
&  &  &  &  &  & A_{2}^{\left(  2\right)  } & A_{1}^{\left(  2\right)  } &
A_{0}^{\left(  2\right)  } &  & \\
&  &  &  &  &  &  & A_{2}^{\left(  3\right)  } & A_{1}^{\left(  3\right)  } &
A_{0}^{\left(  3\right)  } & \\
&  &  &  &  &  &  &  & \ddots & \ddots & \ddots
\end{array}
\right)  ,
\]
where $\widetilde{A}_{0}^{\left(  0\right)  }$, $\widetilde{C}$,
$\widetilde{A}_{2}^{\left(  1\right)  }$, $\widetilde{B}_{0}^{\left(
0\right)  }$ and $\widetilde{B}_{i}^{\left(  -l\right)  }$ for $i=0,1,2$ and
$l\leq-1$, can be written easily and their details are omitted here.

The following theorem shows that the first passage time $\xi$ is of
bidirectional infinite phase type, and thus provides a method to compute the
average first passage time $E\left[  \xi\right]  $. The proof is easy and is
omitted here.

\begin{The}
\label{The3}If the arriving A-customer enters State $\left(  i-1,j\right)  $
with probability $\phi_{i,j}$ for $i\geq1$ at time $0$, then the first passage
time $\xi$ is of bidirectional infinite phase type with an irreducible
representation $\left(  \Phi,T\right)  $, and%
\begin{equation}
E\left[  \xi\right]  =-\Phi T^{-1}e\text{,} \label{S-8}%
\end{equation}
where $\Phi$ is given in (\ref{S-13}), $T^{-1}$ is the maximal non-positive
inverse of the matrix $T$.
\end{The}

When the PH distribution is unilateral infinite, Chapter 8 of Li
\cite{Li:2010} provides a detailed discussion for how to compute the maximal
non-positive inverse matrix $T^{-1}$ of infinite sizes by means of the
RG-factorizations given in Li \cite{Li:2010}. Furthermore, if $T$ is the
infinitesmall generator of an irreducible QBD process, then the maximal
non-positive inverse matrix $T^{-1}$ can be explicitly expressed by means of
the $R$-, $U$- and $G$-measures$.$

When the PH distribution is bidirectional infinite, to our best knowledge,
this paper is the first to deal with the maximal non-positive inverse matrix
$T^{-1}$ of bidirectional infinite sizes. To this end, we write%
\[
T=\left(
\begin{array}
[c]{cc}%
T_{1,1} & T_{1,2}\\
T_{2,1} & T_{2,2}%
\end{array}
\right)  ,
\]
where%
\[
T_{1,1}=\left(
\begin{array}
[c]{ccccc}%
\ddots & \ddots & \ddots &  & \\
& \widetilde{B}_{0}^{\left(  -2\right)  } & \widetilde{B}_{1}^{\left(
-2\right)  } & \widetilde{B}_{2}^{\left(  -2\right)  } & \\
&  & \widetilde{B}_{0}^{\left(  -1\right)  } & \widetilde{B}_{1}^{\left(
-1\right)  } & \widetilde{B}_{2}^{\left(  -1\right)  }\\
&  &  & \widetilde{B}_{0}^{\left(  0\right)  } & \widetilde{C}%
\end{array}
\right)  ,\text{ }T_{1,2}=\left(
\begin{array}
[c]{ccccc}
&  &  &  & \\
&  &  &  & \\
&  &  &  & \\
\widetilde{A}_{0}^{\left(  0\right)  } &  &  &  &
\end{array}
\right)  ,
\]%
\[
T_{2,1}=\left(
\begin{array}
[c]{ccccc}
&  &  &  & \widetilde{A}_{2}^{\left(  1\right)  }\\
&  &  &  & \\
&  &  &  & \\
&  &  &  &
\end{array}
\right)  ,\text{ }T_{2,2}=\left(
\begin{array}
[c]{ccccc}%
A_{1}^{\left(  1\right)  } & A_{0}^{\left(  1\right)  } &  &  & \\
A_{2}^{\left(  2\right)  } & A_{1}^{\left(  2\right)  } & A_{0}^{\left(
2\right)  } &  & \\
& A_{2}^{\left(  3\right)  } & A_{1}^{\left(  3\right)  } & A_{0}^{\left(
3\right)  } & \\
&  & \ddots & \ddots & \ddots
\end{array}
\right)  .
\]

Note that the Markov chain $Q$ is irreducible and $T_{2,1}e\gneq0$ due to
$\widetilde{A}_{2}^{\left(  1\right)  }e\gneq0$, thus the submatrix $T_{2,2}$
must be invertible. Based on this, it is easy to check that%

\begin{equation}
T^{-1}=\left[
\begin{array}
[c]{c}%
T_{1,1;2}^{-1}\text{ \ \ \ \ \ \ \ \ \ }-T_{1,1;2}^{-1}T_{1,2}T_{2,2}^{-1}\\
-T_{2,2}^{-1}T_{2,1}T_{1,1;2}^{-1}\text{ \ \ \ \ }T_{2,2}^{-1}+T_{2,2}%
^{-1}T_{2,1}T_{1,1;2}^{-1}T_{1,2}T_{2,2}^{-1}%
\end{array}
\right]  , \label{S-10}%
\end{equation}
where%
\begin{equation}
T_{1,1;2}=T_{1,1}-T_{1,2}T_{2,2}^{-1}T_{2,1}. \label{S-11}%
\end{equation}
It is easy to see that the maximal non-positive inverse matrix $T^{-1}$ of
bidirectional infinite sizes can be expressed by means of the maximal
non-positive inverse matrix $T_{2,2}^{-1}$ of unilateral infinite sizes. This
relation plays a key role in setting up the PH distribution of bidirectional
infinite sizes.

In what follows, we simply discuss the maximal non-positive inverse matrix
$T_{2,2}^{-1}$ of unilateral infinite sizes by means of the RG-factorizations,
e.g., see Chapter 1 of Li \cite{Li:2010} or Li and Cao \cite{Li:2004} for more details.

For the QBD process $T_{2,2}$ of unilateral infinite sizes, we define the
UL-type $\mathbf{U}$-, $\mathbf{R}$- and $\mathbf{G}$-measures as%
\[
\mathbf{U}_{k}=A_{1}^{\left(  k\right)  }+A_{0}^{\left(  k\right)  }\left(
-\mathbf{U}_{k+1}^{-1}\right)  A_{2}^{\left(  k+1\right)  },\text{ \ }k\geq1,
\]%

\[
\mathbf{R}_{k}=A_{0}^{\left(  k\right)  }\left(  -\mathbf{U}_{k+1}%
^{-1}\right)  ,\text{ \ }k\geq1,
\]
and%
\[
\mathbf{G}_{k}=\left(  -\mathbf{U}_{k}^{-1}\right)  A_{2}^{\left(  k\right)
},\text{ \ }k\geq2.
\]
Note that the matrix sequence $\left\{  \mathbf{R}_{k},k\geq1\right\}  $ is
the minimal nonnegative solution to the system of nonlinear matrix equations%
\[
A_{0}^{\left(  k\right)  }+\mathbf{R}_{k}A_{1}^{\left(  k+1\right)
}+\mathbf{R}_{k}\mathbf{R}_{k+1}A_{2}^{\left(  k+2\right)  }=0,\text{ \ }%
k\geq1;
\]
and the matrix sequence $\left\{  \mathbf{G}_{k},k\geq2\right\}  $ is the
minimal nonnegative solution to the system of nonlinear matrix equations%
\[
A_{0}^{\left(  k\right)  }\mathbf{G}_{k+1}\mathbf{G}_{k}+A_{1}^{\left(
k\right)  }\mathbf{G}_{k}+A_{2}^{\left(  k\right)  }=0,\text{ \ }k\geq2.
\]
Based on this, the UL-type RG-factorization of the matrix $T_{2,2}$ of
unilateral infinite sizes is given by%
\begin{equation}
T_{2,2}=\left(  I-\mathbf{R}_{U}\right)  \mathbf{U}_{D}\left(  I-\mathbf{G}%
_{L}\right)  , \label{S-9}%
\end{equation}
where%

\[
I-\mathbf{R}_{U}=\left(
\begin{array}
[c]{ccccc}%
I & -\mathbf{R}_{1} &  &  & \\
& I & -\mathbf{R}_{2} &  & \\
&  & I & -\mathbf{R}_{3} & \\
&  &  & I & \ddots\\
&  &  &  & \ddots
\end{array}
\right)  ,
\]%
\[
\mathbf{U}_{D}=\text{diag}\left(  \mathbf{U}_{1},\mathbf{U}_{2},\mathbf{U}%
_{3},\ldots\right)  ,
\]%
\[
I-\mathbf{G}_{L}=\left(
\begin{array}
[c]{ccccc}%
I &  &  &  & \\
-\mathbf{G}_{2} & I &  &  & \\
& -\mathbf{G}_{3} & I &  & \\
&  & -\mathbf{G}_{4} & I & \\
&  &  & \ddots & \ddots
\end{array}
\right)  .
\]

By using $\mathbf{R}$-measure $\left\{  \mathbf{R}_{k}:k\geq1\right\}  $ and
the $\mathbf{G}$-measure $\left\{  \mathbf{G}_{k}:k\geq2\right\}  $, we write%
\[
X_{k}^{\left(  k\right)  }=I,
\]%
\[
X_{k+l}^{\left(  k\right)  }=\mathbf{R}_{k}\mathbf{R}_{k+1}\cdots
\mathbf{R}_{k+l-1},\text{ }k\geq1,\text{ }l\geq1,
\]%
\[
Y_{k}^{\left(  k\right)  }=I,
\]%
\[
Y_{k-l}^{\left(  k\right)  }=\mathbf{G}_{k}\mathbf{G}_{k-1}\cdots
\mathbf{G}_{k-l+1},\text{ }k>l\geq1.
\]
From the UL-type RG-factorization (\ref{S-9}), it is easy to check that the
three matrices $I-\mathbf{R}_{U}$, $\mathbf{U}_{D}$ and $I-\mathbf{G}_{L}$ are
all invertible, and%
\[
\left(  I-\mathbf{R}_{U}\right)  ^{-1}=\left(
\begin{array}
[c]{ccccc}%
X_{1}^{\left(  1\right)  } & X_{2}^{\left(  1\right)  } & X_{3}^{\left(
1\right)  } & X_{4}^{\left(  1\right)  } & \ldots\\
& X_{2}^{\left(  2\right)  } & X_{3}^{\left(  2\right)  } & X_{4}^{\left(
2\right)  } & \ldots\\
&  & X_{3}^{\left(  3\right)  } & X_{4}^{\left(  3\right)  } & \ldots\\
&  &  & X_{4}^{\left(  4\right)  } & \ldots\\
&  &  &  & \ddots
\end{array}
\right)  ,
\]%
\[
\mathbf{U}_{D}^{-1}=\text{diag}\left(  \mathbf{U}_{1}^{-1},\mathbf{U}_{2}%
^{-1},\mathbf{U}_{3}^{-1},\mathbf{U}_{4}^{-1},\ldots\right)  ,
\]%
\[
\left(  I-\mathbf{G}_{L}\right)  ^{-1}=\left(
\begin{array}
[c]{ccccc}%
Y_{1}^{\left(  1\right)  } &  &  &  & \\
Y_{1}^{\left(  2\right)  } & Y_{2}^{\left(  2\right)  } &  &  & \\
Y_{1}^{\left(  3\right)  } & Y_{2}^{\left(  3\right)  } & Y_{3}^{\left(
3\right)  } &  & \\
Y_{1}^{\left(  4\right)  } & Y_{2}^{\left(  4\right)  } & Y_{3}^{\left(
4\right)  } & Y_{4}^{\left(  4\right)  } & \\
\vdots & \vdots & \vdots & \vdots & \ddots
\end{array}
\right)  .
\]
We obtain%
\begin{equation}
T_{2,2}^{-1}=\left(  I-\mathbf{G}_{L}\right)  ^{-1}\mathbf{U}_{D}^{-1}\left(
I-\mathbf{R}_{U}\right)  ^{-1}. \label{S-12}%
\end{equation}

By using (\ref{S-10}), (\ref{S-11}) and (\ref{S-12}), we obtain the maximal
non-positive inverse matrix $T^{-1}$ of bidirectional infinite sizes.
Therefore, we can compute the average first passage time $E\left[  \xi\right]
$, and $E\left[  W\right]  \leq E\left[  \xi\right]  $.

\section{The Departure Process}

In this section, we use the MMAP of bidirectional infinite sizes to discuss
the departure process with three types of customers (i.e., impatient
A-customers, impatient B-customers, and groups of $m$ A-customers and $n$
B-customers) in the matched queue with matching batch pair $(m,n)$.

To set up the MMAP of bidirectional infinite sizes, we need to establish four
matrices: $D_{0}$, $D_{A}$, $D_{B}$ and $D_{AB}$, where $D_{A}$ and $D_{B}$
are the departure-rate matrices of A-customers and B-customers due to their
impatient behavior, respectively; $D_{AB}$ is the departure-rate matrix of
groups matched by $m$ A-customers and $n$ B-customers; and $D_{0}$ is the
transition rate matrix of stochastic environment in the MMAP. From the
infinitesimal generator (\ref{Q}), we have%
\begin{equation}
D_{0}=Q-D_{A}-D_{B}-D_{AB}, \label{D0}%
\end{equation}
where%

\begin{equation}
D_{A}=\left(
\begin{array}
[c]{ccccccc}%
\ddots &  &  &  &  &  & \\
\ddots & D_{A}^{\left(  -2\right)  } &  &  &  &  & \\
& 0 & D_{A}^{\left(  -1\right)  } &  &  &  & \\
&  & 0 & D_{A}^{\left(  0\right)  } &  &  & \\
&  &  & \widetilde{D}_{A}^{\left(  1\right)  } & D_{A}^{\left(  1\right)  } &
& \\
&  &  &  & \widetilde{D}_{A}^{\left(  2\right)  } & D_{A}^{\left(  2\right)  }
& \\
&  &  &  &  & \ddots & \ddots
\end{array}
\right)  , \label{DA}%
\end{equation}%
\[
D_{A}^{\left(  k\right)  }=\left(
\begin{array}
[c]{ccccc}%
0 &  &  &  & \\
A_{1,1}^{\left(  k\right)  } & 0 &  &  & \\
& \ddots & \ddots &  & \\
&  & A_{1,m-2}^{\left(  k\right)  } & 0 & \\
&  &  & A_{1,m-1}^{\left(  k\right)  } & 0
\end{array}
\right)  _{mn\times mn},k\geq0,
\]%
\[
A_{1,i}^{\left(  k\right)  }=\left(
\begin{array}
[c]{cccc}%
\left(  km+i\right)  \theta_{1} &  &  & \\
& \left(  km+i\right)  \theta_{1} &  & \\
&  & \ddots & \\
&  &  & \left(  km+i\right)  \theta_{1}%
\end{array}
\right)  _{n\times n},1\leq i\leq m-1,
\]%

\[
D_{A}^{\left(  l\right)  }=\left(
\begin{array}
[c]{ccccc}%
A_{1,1} &  &  &  & \\
& A_{1,2} &  &  & \\
&  & \ddots &  & \\
&  &  & A_{1,n-1} & \\
&  &  &  & A_{1,n}%
\end{array}
\right)  _{nm\times nm},\text{ }l\leq-1,
\]%
\[
A_{1,i}=\left(
\begin{array}
[c]{ccccc}%
0 & \left(  m-1\right)  \theta_{1} &  &  & \\
& 0 & \left(  m-2\right)  \theta_{1} &  & \\
&  & 0 & \ddots & \\
&  &  & \ddots & \theta_{1}\\
&  &  &  & 0
\end{array}
\right)  _{m\times m},\text{ }1\leq i\leq n,
\]%
\[
\widetilde{D}_{A}^{\left(  k\right)  }=\left(
\begin{array}
[c]{cc}%
{\Large 0}_{n\times\left[  \left(  m-1\right)  n\right]  } & \left.
\begin{array}
[c]{ccc}%
km\theta_{1} &  & \\
& \ddots & \\
&  & km\theta_{1}%
\end{array}
\right\}  n\times n\\
& {\Large 0}_{\left[  \left(  m-1\right)  n\right]  \times n}\text{ }%
\end{array}
\right)  _{mn\times mn},\text{\ }k\geq1;
\]%

\begin{equation}
D_{B}=\left(
\begin{array}
[c]{ccccccc}%
\ddots & \ddots &  &  &  &  & \\
& D_{B}^{\left(  -2\right)  } & \widetilde{D}_{B}^{\left(  -2\right)  } &  &
&  & \\
&  & D_{B}^{\left(  -1\right)  } & \widetilde{D}_{B}^{\left(  -1\right)  } &
&  & \\
&  &  & D_{B}^{\left(  0\right)  } & 0 &  & \\
&  &  &  & D_{B}^{\left(  1\right)  } & 0 & \\
&  &  &  &  & D_{B}^{\left(  2\right)  } & \ddots\\
&  &  &  &  &  & \ddots
\end{array}
\right)  , \label{DB}%
\end{equation}%
\[
D_{B}^{\left(  k\right)  }=\left(
\begin{array}
[c]{ccccc}%
B_{1,1} &  &  &  & \\
& B_{1,2} &  &  & \\
&  & \ddots &  & \\
&  &  & B_{1,m-1} & \\
&  &  &  & B_{1,m}%
\end{array}
\right)  _{mn\times mn},k\geq0,
\]%
\[
B_{1,i}=\left(
\begin{array}
[c]{ccccc}%
0 &  &  &  & \\
\theta_{2} & \ddots &  &  & \\
& \ddots & 0 &  & \\
&  & \left(  n-2\right)  \theta_{2} & 0 & \\
&  &  & \left(  n-1\right)  \theta_{2} & 0
\end{array}
\right)  _{n\times n},\text{ }1\leq i\leq m,
\]%

\[
D_{B}^{\left(  l\right)  }=\left(
\begin{array}
[c]{ccccc}%
0 & B_{1,1}^{\left(  l\right)  } &  &  & \\
& 0 & B_{1,2}^{\left(  l\right)  } &  & \\
&  & \ddots & \ddots & \\
&  &  & 0 & B_{1,n-1}^{\left(  l\right)  }\\
&  &  &  & 0
\end{array}
\right)  _{nm\times nm},\text{ }l\leq-1,
\]%
\[
B_{1,i}^{\left(  l\right)  }=\left(
\begin{array}
[c]{cccc}%
\left(  -ln+n-i\right)  \theta_{2} &  &  & \\
& \left(  -ln+n-i\right)  \theta_{2} &  & \\
&  & \ddots & \\
&  &  & \left(  -ln+n-i\right)  \theta_{2}%
\end{array}
\right)  _{m\times m},1\leq i\leq n-1,
\]%
\[
\widetilde{D}_{B}^{\left(  -1\right)  }=\left(
\begin{array}
[c]{cccc}%
{\Large 0}_{\left[  \left(  n-1\right)  m\right]  \times n} & {\Large 0}%
_{\left[  \left(  n-1\right)  m\right]  \times n} & \cdots & {\Large 0}%
_{\left[  \left(  n-1\right)  m\right]  \times n}\\
\text{\ }\left.
\begin{array}
[c]{cc}%
0 & 0\\
\vdots & \vdots\\
0 & 0\\
0 & n\theta_{2}%
\end{array}
\right\}  m\times n & \text{\ }\left.
\begin{array}
[c]{cc}%
0 & 0\\
\vdots & \vdots\\
0 & n\theta_{2}\\
0 & 0
\end{array}
\right\}  m\times n & \cdots & \text{ \ }\left.
\begin{array}
[c]{cc}%
0 & n\theta_{2}\\
\vdots & \vdots\\
0 & 0\\
0 & 0
\end{array}
\right\}  m\times n
\end{array}
\right)  _{mn\times mn},
\]%
\[
\widetilde{D}_{B}^{\left(  l\right)  }=\left(
\begin{array}
[c]{cc}%
{\Large 0}_{\left[  \left(  n-1\right)  m\right]  \times m}\text{\ \ \ \ \ \ }
& {\Large 0}_{\left[  \left(  n-1\right)  m\right]  \times\left[  \left(
n-1\right)  m\right]  }\\
\left.
\begin{array}
[c]{ccc}%
-ln\theta_{2} &  & \\
& \ddots & \\
&  & -ln\theta_{2}%
\end{array}
\right\}  m\times m & {\Large 0}_{m\times\left[  \left(  n-1\right)  m\right]
}%
\end{array}
\right)  _{mn\times mn},\text{ }l\leq-2;
\]%

\begin{equation}
D_{AB}=\left(
\begin{array}
[c]{ccccccccc}%
\ddots & \ddots &  &  &  &  &  &  & \\
\ddots & 0 & D_{AB}^{\left(  -3\right)  } &  &  &  &  &  & \\
& 0 & 0 & D_{AB}^{\left(  -2\right)  } &  &  &  &  & \\
&  & 0 & 0 & D_{AB}^{\left(  -1\right)  } &  &  &  & \\
&  &  & 0 & 0 & 0 &  &  & \\
&  &  &  & D_{AB}^{\left(  1\right)  } & 0 & 0 &  & \\
&  &  &  &  & D_{AB}^{\left(  2\right)  } & 0 & 0 & \\
&  &  &  &  &  & D_{AB}^{\left(  3\right)  } & 0 & \ddots\\
&  &  &  &  &  &  & \ddots & \ddots
\end{array}
\right)  , \label{DAB}%
\end{equation}%
\[
D_{AB}^{\left(  k\right)  }=\left(
\begin{array}
[c]{ccc}%
\left.
\begin{array}
[c]{ccc}%
0 &  & \\
& \ddots & \\
\lambda_{2} &  & 0
\end{array}
\right\}  n\times n &  & \\
& \ddots & \\
&  & \left.
\begin{array}
[c]{ccc}%
0 &  & \\
& \ddots & \\
\lambda_{2} &  & 0
\end{array}
\right\}  n\times n
\end{array}
\right)  _{mn\times mn},\text{\ }k\geq1,
\]%
\[
D_{AB}^{\left(  -1\right)  }=\left(
\begin{array}
[c]{cc}%
\left.
\begin{array}
[c]{cccc}%
0 & \cdots & 0 & \lambda_{1}\\
\vdots &  & \vdots & \vdots\\
0 & \cdots & 0 & 0
\end{array}
\right\}  m\times n & {\Large 0}_{m\times\left[  \left(  m-1\right)  n\right]
}\\
\left.
\begin{array}
[c]{cccc}%
0 & \cdots & \lambda_{1} & 0\\
\vdots &  & \vdots & \vdots\\
0 & \cdots & 0 & 0
\end{array}
\right\}  m\times n & {\Large 0}_{m\times\left[  \left(  m-1\right)  n\right]
}\\
\vdots\text{\ \ } & \vdots\\
\left.
\begin{array}
[c]{cccc}%
\lambda_{1} & \cdots & 0 & 0\\
\vdots &  & \vdots & \vdots\\
0 & \cdots & 0 & 0
\end{array}
\right\}  m\times n & {\Large 0}_{m\times\left[  \left(  m-1\right)  n\right]
}%
\end{array}
\right)  _{mn\times mn},
\]%

\[
D_{AB}^{\left(  l\right)  }=\left(
\begin{array}
[c]{ccc}%
\left.
\begin{array}
[c]{ccc}%
0 &  & \lambda_{1}\\
& \ddots & \\
&  & 0
\end{array}
\right\}  m\times m &  & \\
& \ddots & \\
&  & \left.
\begin{array}
[c]{ccc}%
0 &  & \lambda_{1}\\
& \ddots & \\
&  & 0
\end{array}
\right\}  m\times m
\end{array}
\right)  _{mn\times mn},\text{ }l\leq-2,
\]

The following theorem shows that the departure process of the matched queue
with matching batch pair $(m,n)$ is an MMAP of bidirectional infinite sizes
whose matrix sequence is given by $\left\{  D_{0},D_{A},D_{B},D_{AB}\right\}
$. The proof is easy and omitted here by means of He and Neuts \cite{He:1996,
He:1998}.

\begin{The}
In the matched queue with matching batch pair $(m,n)$, its departure process
has three types of customers: The impatient A- and B-customers, and the groups
matched by $m$ A-customers and $n$ B-customers. Also, the departure process
with the three types of customers is an MMAP of bidirectional infinite sizes
with the matrix sequence $\left\{  D_{0},D_{A},D_{B},D_{AB}\right\}  $.
\end{The}

Note that the matched queue with matching batch pair $(m,n)$ must be stable.
When the stationary probability vector $\pi$ of the Markov process $Q$ is
given, we have

(a) the stationary impatient-departure rate of the A-customers is given by%
\[
\mu_{A}=\pi D_{A}e;
\]

(b) the stationary impatient-departure rate of the B-customers is given by%
\[
\mu_{B}=\pi D_{B}e;
\]

(c) the stationary matched-departure rate of the groups matched by $m$
A-customers and $n$ B-customers is given by%
\[
\mu_{AB}=\pi D_{AB}e.
\]
Also, when observing the A- and B-customers, we have

(d) the stationary departure rate of the A-customers is given by%
\[
\mu_{A}=\pi\left(  D_{A}+mD_{AB}\right)  e;
\]

(e) the stationary departure rate of the B-customers is given by%
\[
\mu_{A}=\pi\left(  D_{B}+nD_{AB}\right)  e;
\]

(f) the stationary departure rate of the A- and B-customers is given by%
\[
\mu_{\text{all}}=\pi\left[  D_{A}+D_{B}+\left(  m+n\right)  D_{AB}\right]  e.
\]

To further understand the departure process with three types of customers,
from the infinitesimal generator (\ref{Q}) and the departure-rate matrices
(\ref{DA}), (\ref{DB}) and (\ref{DAB}), by using (\ref{D0}) we write
\[
D_{0}=\left(
\begin{array}
[c]{ccccccccc}%
\ddots & \ddots & \ddots &  &  &  &  &  & \\
& H_{0}^{\left(  -2\right)  } & H_{1}^{\left(  -2\right)  } & H_{2}^{\left(
-2\right)  } &  &  &  &  & \\
&  & H_{0}^{\left(  -1\right)  } & H_{1}^{\left(  -1\right)  } &
\fbox{$H_{2}^{\left(  -1\right)  }$} &  &  &  & \\
&  &  & \fbox{$H_{0}^{\left(  0\right)  }$} & F & \fbox{$G_{0}^{\left(
0\right)  }$} &  &  & \\
&  &  &  & \fbox{$G_{2}^{\left(  1\right)  }$} & G_{1}^{\left(  1\right)  } &
G_{0}^{\left(  1\right)  } &  & \\
&  &  &  &  & G_{2}^{\left(  2\right)  } & G_{1}^{\left(  2\right)  } &
G_{0}^{\left(  2\right)  } & \\
&  &  &  &  &  & \ddots & \ddots & \ddots
\end{array}
\right)  .
\]
Note that $D_{0}=Q-D_{A}-D_{B}-D_{AB}$ with $D_{A}+D_{B}+D_{AB}\gneq0$, and
the Markov process $Q$ is irreducible and positive recurrent, thus the Markov
process $D_{0}$ is irreducible and transient, so the matrix $D_{0}$ of
bidirectional infinite sizes is invertible. To compute the maximal
non-positive inverse matrix of $D_{0}$, we write%
\[
D_{0}=\left(
\begin{array}
[c]{cc}%
W_{1,1} & W_{1,2}\\
W_{2,1} & W_{2,2}%
\end{array}
\right)  ,
\]
where%
\[
W_{1,1}=\left(
\begin{array}
[c]{ccccc}%
\ddots & \ddots & \ddots &  & \\
& H_{0}^{\left(  -2\right)  } & H_{1}^{\left(  -2\right)  } & H_{2}^{\left(
-2\right)  } & \\
&  & H_{0}^{\left(  -1\right)  } & H_{1}^{\left(  -1\right)  } &
H_{2}^{\left(  -1\right)  }\\
&  &  & H_{0}^{\left(  0\right)  } & F
\end{array}
\right)  ,\text{ }W_{1,2}=\left(
\begin{array}
[c]{ccccc}
&  &  &  & \\
&  &  &  & \\
&  &  &  & \\
G_{0}^{\left(  0\right)  } &  &  &  &
\end{array}
\right)  ,
\]%
\[
W_{2,1}=\left(
\begin{array}
[c]{ccccc}
&  &  &  & G_{2}^{\left(  1\right)  }\\
&  &  &  & \\
&  &  &  & \\
&  &  &  &
\end{array}
\right)  ,\text{ }W_{2,2}=\left(
\begin{array}
[c]{ccccc}%
G_{1}^{\left(  1\right)  } & G_{0}^{\left(  1\right)  } &  &  & \\
G_{2}^{\left(  2\right)  } & G_{1}^{\left(  2\right)  } & G_{0}^{\left(
2\right)  } &  & \\
& G_{2}^{\left(  3\right)  } & G_{1}^{\left(  3\right)  } & G_{0}^{\left(
3\right)  } & \\
&  & \ddots & \ddots & \ddots
\end{array}
\right)  .
\]
Since the Markov chain $Q$ is irreducible and $W_{2,1}e\gneq0$ due to
$G_{2}^{\left(  1\right)  }e\gneq0$, the submatrix $W_{2,2}$ must be
invertible. Based on this, it is easy to check that%

\[
D_{0}^{-1}=\left(
\begin{array}
[c]{c}%
W_{1,1;2}^{-1}\text{ \ \ \ \ \ \ \ \ \ }-W_{1,1;2}^{-1}W_{1,2}W_{2,2}^{-1}\\
-W_{2,2}^{-1}W_{2,1}W_{1,1;2}^{-1}\text{ \ \ \ \ }W_{2,2}^{-1}+W_{2,2}%
^{-1}W_{2,1}W_{1,1;2}^{-1}W_{1,2}W_{2,2}^{-1}%
\end{array}
\right)  ,
\]
where%
\[
W_{1,1;2}=W_{1,1}-W_{1,2}W_{2,2}^{-1}W_{2,1},
\]
and the maximal non-positive inverse matrix $W_{2,2}^{-1}$ of unilateral
infinite sizes can be computed by means of the RG-factorizations, as seen in
(\ref{S-9}).

Let $A$ and $B$ be the impatient departures of the A- and A-customers,
respectively; and $AB$ be a departure of the groups matched by $m$ A-customers
and $n$ B-customers.

The following theorem describes some probability characteristics of the
departure process with the three types of customers, while its proof is easy
by using He and Neuts \cite{He:1996}, and it is omitted here.

\begin{The}
In the matched queue with matching batch pair $(m,n)$, its departure process
with the three types of customers: The MMAP of bidirectional infinite sizes
with the matrix sequence $\left\{  D_{0},D_{A},D_{B},D_{AB}\right\}  $, has
the probability characteristics as follows:

(a) Backward looking: The probability that the last departure before an
arbitrary time $t$ is marked by $i\in\left\{  A,B,AB\right\}  $ is given by
$\pi D_{i}\left(  -D_{0}^{-1}\right)  e$.

(b) Forward looking: The probability that the first departure after an
arbitrary time $t$ is marked by $i\in\left\{  A,B,AB\right\}  $ is given by
$\pi\left(  -D_{0}^{-1}\right)  D_{i}e$.

(c) At the departure: The probability that an arbitrary departure is of type
$i$ is given by%
\[
\frac{\pi D_{i}e}{\pi\left[  D_{A}+D_{B}+D_{AB}\right]  e},\text{ \ }%
i\in\left\{  A,B,AB\right\}  .
\]
\end{The}

Finally, we extend the above probability characteristics of the departure
process with the three types of customers from a departure point to multiple
departure points. We write that $i_{1},i_{2},\ldots,i_{k}\in\left\{
A,B,AB\right\}  $.

\begin{The}
(a) Backward looking: Let%
\[
i_{k}\leftarrow i_{k-1}\leftarrow\cdots\leftarrow i_{2}\leftarrow
i_{1}\Longleftarrow t
\]
Then the probability that the $k$ consecutive departure points before an
arbitrary time $t$ are marked by $i_{1},i_{2},\ldots,i_{k}\in\left\{
A,B,AB\right\}  $ is given by%
\[
\pi\left[  D_{i_{1}}\left(  -D_{0}^{-1}\right)  \right]  \left[  D_{i_{2}%
}\left(  -D_{0}^{-1}\right)  \right]  \cdots\left[  D_{i_{k}}\left(
-D_{0}^{-1}\right)  \right]  e.
\]

(b) Forward looking: Let%
\[
t\Longrightarrow i_{1}\rightarrow i_{2}\rightarrow\cdots\rightarrow
i_{k-1}\rightarrow i_{k}%
\]
Then the probability that the $k$ consecutive departure points after an
arbitrary time $t$ are marked by $i_{1},i_{2},\ldots,i_{k}\in\left\{
A,B,AB\right\}  $ is given by%
\[
\pi\left[  \left(  -D_{0}^{-1}\right)  D_{i_{1}}\right]  \left[  \left(
-D_{0}^{-1}\right)  D_{i_{2}}\right]  \cdots\left[  \left(  -D_{0}%
^{-1}\right)  D_{i_{k}}\right]  e.
\]

\begin{Rem}
For the matched queue with matching batch pair $(m,n)$, the analysis of the
departure process with the three types of customers is described as the MMAP
of bidirectional infinite sizes with the matrix sequence $\left\{  D_{0}%
,D_{A},D_{B},D_{AB}\right\}  $; while the finite-size case of which was
studied in He and Neuts \cite{He:1996}. In the study of MMAPs, we make two
useful advances: (1) The maximal non-positive inverse matrix $D_{0}^{-1}$ of
bidirectional infinite sizes can be expressed by means of the maximal
non-positive inverse matrix $W_{2,2}^{-1}$ of unilateral infinite sizes, and
(2) the maximal non-positive inverse matrix $W_{2,2}^{-1}$ of unilateral
infinite sizes can be determined by means of the R-, U- and G-measures, as
seen in how to compute the maximal non-positive inverse matrix $T^{-1}$ of the
previous section. Therefore, our RG-factorization method generalizes the MMAPs
of finite sizes given in He and Neuts \cite{He:1996} to the MMAPs of (either
unilateral or bidirectional) infinite sizes such that the MMAPs can adapt to a
wider range of practical applications.
\end{Rem}

\begin{Rem}
It is necessary and useful to study the departure process of the matched queue
with matching batch pair $(m,n)$ (and more generally, stochastic models). Such
a departure process can be regarded as a new input in a local node of the
system. See He et al. \cite{He:2007} for more details. This is very useful in
the study of large-scale stochastic networks.
\end{Rem}
\end{The}

\section{Concluding Remarks}

In this paper, we discuss an interesting but challenging bilateral
stochastically matching problem: A more general matched queue with matching
batch pair $(m,n)$ and two types of impatient customers. We show that the
matched queue with matching batch pair $(m,n)$ can be expressed as a novel
level-dependent QBD process with bidirectional infinite levels. Based on this,
we provide a detailed analysis for this matched queue, including the system
stability, the average stationary queue lengths, the average sojourn times,
and the departure process. We believe that the methodology and results
developed in this paper can be applicable to analyze more general matched
queues, which are widely encountered in many practical areas, for example,
sharing economy, ridesharing platform, bilateral market, organ
transplantation, taxi services, assembly systems, and so on.

Along these lines, we will continue our future research on the following directions:

-- Considering probabilistic matching system with matching pair $\left(
1,1\right)  $.

-- Studying probabilistic matching system with matching batch pair $\left(
m,n\right)  $.

-- Discussing some Phase-type and MAP factors in the study of matched queues.

-- Analyzing fluid and diffusion approximations for matched system with
matching batch pair $\left(  m,n\right)  $.

-- Developing stochastic optimization, and Markov decision processes in the
study of matched queues.

\section*{Acknowledgements}

Quan-Lin Li was supported by the National Natural Science Foundation of China
under grants No. 71671158 and 71932002 and by Beijing Social Science
Foundation Research Base Project under grant No. 19JDGLA004.

\section*{\textbf{Appendix}}

The appendix contains the proofs of Theorems \ref{The1} and \ref{The2}.

\textbf{Proof of Theorem \ref{The1}.}\textit{ }It is clear that the QBD
process $Q$ with bidirectional infinite levels is irreducible through
observing Figures 1 and 2, since this matched queue contains two Poisson
inputs and two exponential impatient times.

From returning to Level $0$, it is easy to see that the QBD process $Q$ with
bidirectional infinite levels is positive recurrent if and only if the two
unilateral QBD processes $Q_{A}$ and $Q_{B}$ are both positive recurrent.
Thus, our aim is to prove that if $\left(  \theta_{1},\theta_{2}\right)  >0$,
then the two unilateral QBD processes $Q_{A}$ and $Q_{B}$ are both positive
recurrent by means of the mean drift technique by Neuts [48] and Li [42].

For the QBD process $Q_{A}$,\textit{ }let $\mathbb{A}_{k}=A_{0}^{\left(
k\right)  }+A_{1}^{\left(  k\right)  }+A_{2}^{\left(  k\right)  }$ in Level
$k$. Then%
\[
\mathbb{A}_{k}=\left(
\begin{array}
[c]{ccccc}%
A_{1,1,1} & A_{2,1} &  &  & A_{2,2,2}\\
A_{3,1} & A_{1,2,1} & A_{2,2} &  & \\
& \ddots & \ddots & \ddots & \\
&  & A_{3,m-2} & A_{1,m-1,1} & A_{2,m-1}\\
A_{0} &  &  & A_{3,m-1} & A_{1,m,1}%
\end{array}
\right)  _{mn\times mn},
\]
where%
\[
A_{0}=\left(
\begin{array}
[c]{ccc}%
\lambda_{1} &  & \\
& \ddots & \\
&  & \lambda_{1}%
\end{array}
\right)  _{n\times n},\text{ }A_{2,2,2}=\left(
\begin{array}
[c]{ccc}%
km\theta_{1} &  & \\
& \ddots & \\
&  & km\theta_{1}%
\end{array}
\right)  _{n\times n},
\]%
\[
A_{1,i,1}=\left(
\begin{array}
[c]{ccccc}%
a_{1,i}^{\left(  1\right)  } & \lambda_{2} &  &  & \\
\theta_{2} & a_{1,i}^{\left(  2\right)  } & \lambda_{2} &  & \\
& \ddots & \ddots & \ddots & \\
&  & \left(  n-2\right)  \theta_{2} & a_{1,i}^{\left(  n-1\right)  } &
\lambda_{2}\\
\lambda_{2} &  &  & \left(  n-1\right)  \theta_{2} & a_{1,i}^{\left(
n\right)  }%
\end{array}
\right)  _{n\times n},\text{ }1\leq i\leq m,
\]%
\[
a_{1,i}^{\left(  r\right)  }=-\left(  \lambda_{1}+\lambda_{2}+\left(
r-1\right)  \theta_{2}+\left(  km+i-1\right)  \theta_{1}\right)  ,\text{
}1\leq r\leq n,
\]%
\[
A_{2,i}=\left(
\begin{array}
[c]{cccc}%
\lambda_{1} &  &  & \\
& \lambda_{1} &  & \\
&  & \ddots & \\
&  &  & \lambda_{1}%
\end{array}
\right)  _{n\times n},1\leq i\leq m-1,
\]%
\[
A_{3,i}=\left(
\begin{array}
[c]{cccc}%
\left(  km+i\right)  \theta_{1} &  &  & \\
& \left(  km+i\right)  \theta_{1} &  & \\
&  & \ddots & \\
&  &  & \left(  km+i\right)  \theta_{1}%
\end{array}
\right)  _{n\times n},1\leq i\leq m-1.
\]

Let $\alpha=\left(  \alpha_{1},\alpha_{2},...,\alpha_{n};...;\alpha_{\left(
m-1\right)  n+1},\alpha_{\left(  m-1\right)  n+2},...,\alpha_{\left(
m-1\right)  n+n}\right)  $ be the stationary probability vector of the Markov
process $\mathbb{A}_{k}$. Then
\[
\alpha\mathbb{A}_{k}=\mathbf{0},\text{ \ }\alpha\mathbf{e}=1.
\]
Note that $\alpha>0$, since the Markov process $\mathbb{A}_{k}$ is irreducible.

Once the stationary probability vector $\alpha$ is obtained, we can compute
the (upward and downward) mean drift rates of the QBD process $Q_{A}$. From
Level $k$ to Level $k+1$, the upward mean drift rate is given by%
\[
\alpha A_{0}^{\left(  k\right)  }\mathbf{e}=\lambda_{1}\left(  \alpha_{\left(
m-1\right)  n+1}+\alpha_{\left(  m-1\right)  n+2}+\cdots+\alpha_{\left(
m-1\right)  n+n}\right)  .
\]
Similarly, from Level $k$ to Level $k-1$, the downward mean drift rate is
given by%
\[
\text{\ }\alpha A_{2}^{\left(  k\right)  }\mathbf{e}=\lambda_{2}\left(
\alpha_{n}+\alpha_{2n}+\cdots+\alpha_{mn}\right)  +km\theta_{1}\left(
\alpha_{1}+\alpha_{2}+\cdots+\alpha_{n}\right)  .
\]
Note that $k$ is a positive integer, $\lambda_{1}>0$ and $\theta_{1}>0$, it is
easy to check that if $k>\max\left\{  1,\lambda_{1}/m\theta_{1}\right\}  $,
then $\alpha A_{0}^{\left(  k\right)  }\mathbf{e<}$ $\alpha A_{2}^{\left(
k\right)  }\mathbf{e}$.\textbf{ }Therefore, the QBD process $Q_{A}$ is
positive recurrent due to the fact that the mean drift rates: $\alpha
A_{0}^{\left(  k\right)  }\mathbf{e<}$ $\alpha A_{2}^{\left(  k\right)
}\mathbf{e}$ for a bigger positive integer $k$, this can hold because $k$ goes
to infinity.

Similarly, we discuss the stability of the QBD process $Q_{B}$. Let
$\mathbb{B}_{l}=B_{0}^{\left(  l\right)  }+B_{1}^{\left(  l\right)  }%
+B_{2}^{\left(  l\right)  }$ in Level $l$ for $l\leq-2$. Then%
\[
\mathbb{B}_{k}=\left(
\begin{array}
[c]{ccccc}%
B_{1,1,1} & B_{2,1} &  &  & B_{0}\\
B_{3,1} & B_{1,2,1} & B_{2,2} &  & \\
& \ddots & \ddots & \ddots & \\
&  & B_{3,n-2} & B_{1,n-1,1} & B_{2,n-1}\\
B_{2,2,2} &  &  & B_{3,n-1} & B_{1,n,1}%
\end{array}
\right)  _{mn\times mn},\text{ }%
\]
where%
\[
B_{0}=\left(
\begin{array}
[c]{ccc}%
\lambda_{2} &  & \\
& \ddots & \\
&  & \lambda_{2}%
\end{array}
\right)  _{m\times m},\text{ }B_{2,2,2}=\left(
\begin{array}
[c]{ccc}%
-ln\theta_{2} &  & \\
& \ddots & \\
&  & -ln\theta_{2}%
\end{array}
\right)  _{m\times m},
\]%
\[
B_{1,i,1}=\left(
\begin{array}
[c]{ccccc}%
b_{1,i}^{\left(  1\right)  } & \left(  m-1\right)  \theta_{1} &  &  &
\lambda_{1}\\
\lambda_{1} & b_{1,i}^{\left(  2\right)  } & \left(  m-2\right)  \theta_{1} &
& \\
& \ddots & \ddots & \ddots & \\
&  & \lambda_{1} & b_{1,i}^{\left(  m-1\right)  } & \theta_{1}\\
&  &  & \lambda_{1} & b_{1,i}^{\left(  m\right)  }%
\end{array}
\right)  _{m\times m},\text{ }1\leq i\leq n,
\]%
\[
b_{1,i}^{\left(  r\right)  }=-\left(  \lambda_{1}+\lambda_{2}+\left(
m-r\right)  \theta_{1}+\left(  -ln+n-i\right)  \theta_{2}\right)  ,\text{
}1\leq r\leq m,
\]%
\[
B_{2,i}=\left(
\begin{array}
[c]{cccc}%
\left(  -ln+n-i\right)  \theta_{2} &  &  & \\
& \left(  -ln+n-i\right)  \theta_{2} &  & \\
&  & \ddots & \\
&  &  & \left(  -ln+n-i\right)  \theta_{2}%
\end{array}
\right)  _{m\times m},1\leq i\leq n-1,
\]%
\[
B_{3,i}=\left(
\begin{array}
[c]{cccc}%
\lambda_{2} &  &  & \\
& \lambda_{2} &  & \\
&  & \ddots & \\
&  &  & \lambda_{2}%
\end{array}
\right)  _{m\times m},1\leq i\leq n-1.
\]

Let $\beta=\left(  \beta_{1},\beta_{2},...,\beta_{m};...;\beta_{\left(
n-1\right)  m+1},\beta_{\left(  n-1\right)  m+2},...,\beta_{\left(
n-1\right)  m+m}\right)  $ be the stationary probability vector of the Markov
process $\mathbb{B}_{l}$. Then
\[
\beta\mathbb{B}_{l}=\mathbf{0},\text{ \ }\beta\mathbf{e}=1.
\]
Now, we compute the (upward and downward) mean drift rates of the QBD process
$Q_{B}$. From Level $l$ to Level $l-1$, the upward mean drift rate is given by%
\[
\beta B_{0}^{\left(  l\right)  }\mathbf{e}=\lambda_{2}\left(  \beta_{1}%
+\beta_{2}+\cdots+\beta_{m}\right)  .
\]
Similarly, from Level $l$ to Level $l+1$, the downward mean drift rate is
given by%
\begin{align*}
\text{\ }\beta B_{2}^{\left(  l\right)  }\mathbf{e}=  &  \lambda_{1}\left(
\beta_{1}+\beta_{m+1}+\cdots+\beta_{\left(  n-1\right)  m+1}\right) \\
&  -ln\theta_{2}\left(  \beta_{\left(  n-1\right)  m+1}+\beta_{\left(
n-1\right)  m+2}+\cdots+\beta_{\left(  n-1\right)  m+m}\right)  ,
\end{align*}
Note that $l$ is a negative integer, $\lambda_{2}>0$ and $\theta_{2}>0$, it is
easy to check that if $l<-\max\left\{  1,\lambda_{2}/n\theta_{2}\right\}  $,
then $\beta B_{0}^{\left(  l\right)  }\mathbf{e<}$ $\beta B_{2}^{\left(
l\right)  }\mathbf{e}$. This holds because $l$ goes to negative
infinity.\textbf{ }Therefore, the QBD process $Q_{B}$ is positive recurrent.

Based on the above two analysis, the two QBD processes $Q_{A}$ and $Q_{B}$ are
both positive recurrent. Thus the QBD process $Q$ with bidirectional infinite
levels is irreducible and positive recurrent. This further shows that the
matched queue with matching batch pair $(m,n)$ and impatient customers is
stable. This completes the proof. \hfill$\blacksquare$

\textbf{Proof of Theorem \ref{The2}.}\textit{ }The proof is easy through
checking whether $\pi$ satisfies the system of linear equations: $\pi
Q=\mathbf{0}$ and $\pi\mathbf{e}=1$.\ To this end, we consider the following
three different cases:

\textit{Case one:} $k\geq2$. In this case, we need to check that%
\begin{equation}
\pi_{k-1}A_{0}^{\left(  k-1\right)  }+\pi_{k}A_{1}^{\left(  k\right)  }%
+\pi_{k+1}A_{2}^{\left(  k+1\right)  }=0. \label{equ1211}%
\end{equation}
Since $\pi_{k-1}=c\widetilde{\pi}_{k-1}=c\widetilde{\pi}_{1}R_{1}R_{2}\cdots
R_{k-2},\pi_{k}=c\widetilde{\pi}_{k}=c\widetilde{\pi}_{1}R_{1}R_{2}\cdots
R_{k-1}$ and $\pi_{k+1}=c\widetilde{\pi}_{k+1}=c\widetilde{\pi}_{1}R_{1}%
R_{2}\cdots R_{k}$, we obtain%
\begin{align*}
&  \text{ \ \ }\pi_{k-1}A_{0}^{\left(  k-1\right)  }+\pi_{k}A_{1}^{\left(
k\right)  }+\pi_{k+1}A_{2}^{\left(  k+1\right)  }\\
&  =c\widetilde{\pi}_{1}R_{1}R_{2}\cdots R_{k-2}A_{0}^{\left(  k-1\right)
}+c\widetilde{\pi}_{1}R_{1}R_{2}\cdots R_{k-1}A_{1}^{\left(  k\right)
}+c\widetilde{\pi}_{1}R_{1}R_{2}\cdots R_{k}A_{2}^{\left(  k+1\right)  }\\
&  =c\widetilde{\pi}_{1}R_{1}R_{2}\cdots R_{k-2}\left(  A_{0}^{\left(
k-1\right)  }+R_{k-1}A_{1}^{\left(  k\right)  }+R_{k-1}R_{k}A_{2}^{\left(
k+1\right)  }\right)  =0
\end{align*}
by means of (\ref{equ3}).

\textit{Case two:} $k\leq-2$. In this case, we need to check that
\begin{equation}
\pi_{k+1}B_{0}^{\left(  k+1\right)  }+\pi_{k}B_{1}^{\left(  k\right)  }%
+\pi_{k-1}B_{2}^{\left(  k-1\right)  }=0. \label{equ1212}%
\end{equation}
Note that $\pi_{k+1}=c\widetilde{\pi}_{k+1}=c\widetilde{\pi}_{-1}%
\mathbb{R}_{-1}\mathbb{R}_{-2}\cdots\mathbb{R}_{k+2},\pi_{k}=c\widetilde{\pi
}_{k}=c\widetilde{\pi}_{-1}\mathbb{R}_{-1}\mathbb{R}_{-2}\cdots\mathbb{R}%
_{k+1}$ and $\pi_{k-1}=c\widetilde{\pi}_{k-1}=c\widetilde{\pi}_{-1}%
\mathbb{R}_{-1}\mathbb{R}_{-2}\cdots\mathbb{R}_{k}$, we have%
\begin{align*}
&  \text{ \ \ }\pi_{k+1}B_{0}^{\left(  k+1\right)  }+\pi_{k}B_{1}^{\left(
k\right)  }+\pi_{k-1}B_{2}^{\left(  k-1\right)  }\\
&  =c\widetilde{\pi}_{-1}\mathbb{R}_{-1}\mathbb{R}_{-2}\cdots\mathbb{R}%
_{k+2}B_{0}^{\left(  k+1\right)  }+c\widetilde{\pi}_{-1}\mathbb{R}%
_{-1}\mathbb{R}_{-2}\cdots\mathbb{R}_{k+1}B_{1}^{\left(  k\right)
}+c\widetilde{\pi}_{-1}\mathbb{R}_{-1}\mathbb{R}_{-2}\cdots\mathbb{R}_{k}%
B_{2}^{\left(  k-1\right)  }\\
&  =c\widetilde{\pi}_{-1}\mathbb{R}_{-1}\mathbb{R}_{-2}\cdots\mathbb{R}%
_{k+2}\left(  B_{0}^{\left(  k+1\right)  }+\mathbb{R}_{k+1}B_{1}^{\left(
k\right)  }+\mathbb{R}_{k+1}\mathbb{R}_{k}B_{2}^{\left(  k-1\right)  }\right)
=0
\end{align*}
in terms of (\ref{equ7}).

\textit{Case three:} $k=1,$ $0,$ $-1.$ In this case, we obtain%
\begin{equation}
\left\{
\begin{array}
[c]{l}%
\pi_{0}A_{0}^{\left(  0\right)  }+\pi_{1}A_{1}^{\left(  1\right)  }+\pi
_{2}A_{2}^{\left(  2\right)  }=0,\\
\pi_{-1}B_{2}^{\left(  -1\right)  }+\pi_{0}C+\pi_{1}A_{2}^{\left(  1\right)
}=0,\\
\pi_{-2}B_{2}^{\left(  -2\right)  }+\pi_{-1}B_{1}^{\left(  -1\right)  }%
+\pi_{0}B_{0}^{\left(  0\right)  }=0.
\end{array}
\right.  \label{equ1213}%
\end{equation}
Note that $\pi_{0}=c\widetilde{\pi}_{0},\pi_{1}=c\widetilde{\pi}_{1},\pi
_{2}=c\widetilde{\pi}_{2}=c\widetilde{\pi}_{1}R_{1},\pi_{-1}=c\widetilde{\pi
}_{-1}$ and $\pi_{-2}=c\widetilde{\pi}_{-2}=c\widetilde{\pi}_{-1}%
\mathbb{R}_{-1}$, while $\widetilde{\pi}_{-1}$, $\widetilde{\pi}_{0}$ and
$\widetilde{\pi}_{1}$ are given in (\ref{equ12}).

Note that $\pi\mathbf{e}=1$, we have%
\begin{equation}
\sum_{-\infty<k<\infty}\pi_{k}\mathbf{e}=1 \label{equ1214}%
\end{equation}
by means of $\pi_{k}=c\widetilde{\pi}_{k}=c\widetilde{\pi}_{1}R_{1}R_{2}\cdots
R_{k-1}$ for $k\geq2,\pi_{1}=c\widetilde{\pi}_{1},\pi_{0}=c\widetilde{\pi}%
_{0},\pi_{-1}=c\widetilde{\pi}_{-1}$ and $\pi_{k}=c\widetilde{\pi}%
_{k}=c\widetilde{\pi}_{-1}\mathbb{R}_{-1}\mathbb{R}_{-2}\cdots\mathbb{R}%
_{l+1}$ for $l\leq-2$. Thus we have%
\begin{align*}
1  &  =\sum_{-\infty<k<\infty}\pi_{k}\mathbf{e}=\sum_{l\leq-2}\pi
_{l}\mathbf{e+}\pi_{-1}\mathbf{e+}\pi_{0}\mathbf{e+}\pi_{1}\mathbf{e+}%
\sum_{k=2}^{\infty}\pi_{k}\mathbf{e}\\
\text{\ }  &  =\sum_{l\leq-2}c\widetilde{\pi}_{-1}\mathbb{R}_{-1}%
\mathbb{R}_{-2}\cdots\mathbb{R}_{l+1}\mathbf{e+}c\widetilde{\pi}%
_{-1}\mathbf{e+}c\widetilde{\pi}_{0}\mathbf{e+}c\widetilde{\pi}_{1}%
\mathbf{e+}\sum\limits_{k=2}^{\infty}c\widetilde{\pi}_{1}R_{1}R_{2}\cdots
R_{k-1}\mathbf{e}\\
&  =c\left(  \sum_{l\leq-2}\widetilde{\pi}_{-1}\mathbb{R}_{-1}\mathbb{R}%
_{-2}\cdots\mathbb{R}_{l+1}\mathbf{e}+\widetilde{\pi}_{-1}\mathbf{e}%
+\widetilde{\pi}_{0}\mathbf{e}+\widetilde{\pi}_{1}\mathbf{e}+\sum
\limits_{k=2}^{\infty}\widetilde{\pi}_{1}R_{1}R_{2}\cdots R_{k-1}%
\mathbf{e}\right)  ,
\end{align*}
this gives the positive constant $c$ in (\ref{equ121}). This completes the
proof.\hfill$\blacksquare$

\end{document}